
\input amstex
\documentstyle{amsppt}
\magnification=1200
\NoBlackBoxes

\voffset -1cm

\define\al{\alpha}
\define\be{\beta}
\define\ga{\gamma}
\define\Ga{\Gamma}
\define\de{\delta}
\define\la{\lambda}
\define\om{\omega}
\define\Om{\Omega}
\define\ep{\varepsilon}
\define\De{\Delta}

\define\Reg{\operatorname{Reg}}
\define\C{\Bbb C}
\define\Z{\Bbb Z}
\define\R{\Bbb R}
\define\GT{\Bbb{GT}}
\define\DD{\frak D}
\define\HH{\frak H}
\define\Leb{\operatorname{\ell}}
\define\U{\frak U}

\define\m{m}
\define\dual{\widehat{\phantom{a}}}

\define\Dim{\operatorname{Dim}}
\define\diag{\operatorname{diag}}
\define\const{\operatorname{const}}
\define\Herm{\operatorname{Herm}}
\define\wt{\widetilde}

\define\tht{\thetag}

\define\M{\operatorname{Mat}(N,\C)}
\define\Mplus{\M_+}
\define\MM{\M_{>-1}}
\define\Ma{\operatorname{Mat}(N-1,\C)}

\define\MMa{\Ma_{>-1}}
\define\F{\Cal F}
\define\Ex{\operatorname{Ex}}

\define\Dadm{\Cal D_{\operatorname{adm}}}
\define\zw{{z,z,w,w'}}
\define\un{\underline}

\define\Treg{\Cal T_{\operatorname{reg}}}
\define\n{{(N)}}
\define\Zp{\Cal Z_{\operatorname{princ}}}
\define\Zc{\Cal Z_{\operatorname{compl}}}
\define\Zd{\Cal Z_{\operatorname{degen}}}
\define\Zdm{\Cal Z_{\operatorname{degen},m}}
\define\Zdk{\Cal Z_{\operatorname{degen},k}}
\define\Zdl{\Cal Z_{\operatorname{degen},l}}
\define\E{\Bbb E}

\TagsOnRight
\NoRunningHeads

\topmatter

\title The problem of harmonic analysis on the
infinite--dimensional unitary group \endtitle

\author Grigori Olshanski \endauthor

\address Dobrushin Mathematics Laboratory, Institute for Information
Transmission Problems, Bolshoy Karetny 19, 101447 Moscow GSP--4,
Russia 
\endaddress

\email {\tt olsh\@iitp.ru, olsh\@online.ru}
\endemail

\abstract
The goal of harmonic analysis on a (noncommutative) group is to decompose the
most ``natural'' unitary representations of this group (like the regular
representation) on irreducible ones. The infinite--dimensional unitary
group $U(\infty)$ is one of the basic examples of ``big'' groups whose
irreducible representations depend on infinitely many parameters.
Our aim is to explain what the harmonic analysis on $U(\infty)$ consists of.

We deal with unitary representations of a reasonable class, which are
in 1--1 correspondence with characters (central, positive definite,
normalized functions on $U(\infty)$). The decomposition of any
representation of this class is described by a probability 
measure (called spectral measure) on the space of indecomposable
characters. The indecomposable characters were found by Dan Voiculescu
in 1976. 

The main result of the present paper consists in explicitly
constructing a 4--parameter family of ``natural'' representations and
computing their characters. We view these representations as a
substitute of the nonexisting regular representation of $U(\infty)$. 
We state the problem of harmonic analysis on $U(\infty)$ as the
problem of computing the spectral measures for these ``natural''
representations. A solution to this problem is given in the next
paper \cite{BO4}, joint with Alexei Borodin.  

We also prove a few auxiliary general results. In particular, it is
proved that the spectral measure of any character of $U(\infty)$ can
be approximated by a sequence of (discrete) spectral measures for the
restrictions of the character to the compact unitary groups $U(N)$.
This fact is a starting point for computing spectral measures.  

\endabstract

\toc
\widestnumber\head{15}
\head {} Introduction \endhead
\head 1. Characters \endhead
\head 2. Spherical and admissible representations \endhead
\head 3. The representations $T_{z,w}$ \endhead
\head 4. The space $\U$ of virtual unitary matrices and another
realization of the representations $T_{z,w}$ \endhead
\head 5. Measures on the space of Hermitian matrices \endhead
\head 6. The characters $\chi_{z,w}$ \endhead 
\head 7. Other series of characters \endhead
\head 8. Topology on the space of extreme characters \endhead
\head 9. Existence and uniqueness of spectral decomposition \endhead
\head 10. Approximation of spectral measures \endhead
\head 11. Conclusion: the problem of harmonic analysis \endhead
\head {} References \endhead
\endtoc

\endtopmatter

\document

\head Introduction \endhead

\subhead (a) Preface \endsubhead
The problem of noncommutative harmonic analysis consists in
decomposing ``natural'' unitary representations of a given group into
irreducible ones. 

For instance, if $K$ is a compact group, then the
decomposition of the (bi)regular representation in the space $L^2(K)$
is described by the classical Peter--Weyl theorem. 

Another well--known example is the decomposition of the quasiregular
representation of a noncompact simple Lie group $G$ acting in the
$L^2$ space on the Riemannian symmetric space $G/K$.

Note that in the Peter--Weyl theorem, the spectrum of the
decomposition is discrete, while for $L^2(G/K)$ the spectrum is
continuous. The decomposition of $L^2(G/K)$ is described by a
measure living on a region of dimension equal to the rank of the
space $G/K$. 

The problems of decomposing the $L^2$ space on $G$ and on a
pseudo--Riemannian symmetric space $G/H$ belong to the next levels of
difficulty. 

In the present paper and the next one (joint with Alexei Borodin,
\cite{BO4}), we deal with the problem of harmonic
analysis in a totally different situation. The novelty is that the
group is no longer compact or locally compact, its dual space has
infinite dimension, and the decomposition into irreducibles is
governed by a measure with infinite--dimensional support. 

Our main result (established in \cite{BO4}) is an explicit
description of measures which arise in this way (we call them
spectral measures). The description is
given in the language of stochastic point processes. These
probabilistic objects have never emerged in classical representation
theory. However, they happen to be an adequate tool for groups with
infinite--dimensional dual.  

The present paper contains results of two kinds:

First, we construct a family of representations which we consider as
``natural'' ones. In our situation, when the group is not locally
compact, the conventional definition of a regular representation (or
a quasi--regular representation associated with a 
homogeneous space) is not applicable directly. This forces us to choose
another, more sophisticated, way to produce representations. 

Second, we prove necessary general theorems concerning the
spectral measures with infinite--dimensional support. 

This finally allows us to convert the problem of harmonic analysis
to an asymptotic problem of the form which is typical for
random matrix theory or asymptotic combinatorics. We are lead,
however, to a new model, which was not previously examined.

We proceed now to a more detailed description of the contents of the
present paper. 

\subhead (b) The group \endsubhead
Consider the chain of the compact classical groups $U(N)$,
$N=1,2,\dots$, which are embedded one into another in a natural way,
and let $U(\infty)$ be their union. Equivalently, elements of
$U(\infty)$ are infinite unitary matrices $U=[U_{ij}]$, where the
indices $i,j$ take values $1,2,\dots$, and we assume that
$U_{ij}=\delta_{ij}$ for $i+j$ large enough. The group $U(\infty)$ is
one of the fundamental examples of inductive limit groups (another
such example is $S(\infty)$, the union of the finite symmetric
groups). 

Following the philosophy of \cite{Ol1}, \cite{Ol3}, we form a $(G,K)$--pair,
where $G$ is the group $U(\infty)\times U(\infty)$ and $K$ is the
diagonal subgroup in $G$, isomorphic to $U(\infty)$.
This is a Gelfand pair in the sense of \cite{Ol3}.

\subhead (c) Representations and characters \endsubhead
We are dealing with unitary representations $T$ of the group $G$
possessing a distinguished cyclic $K$--invariant vector $\xi$. Such
representations are called {\it spherical representations\/} of the pair
$(G,K)$. They are completely determined by the corresponding
matrix coefficients $\psi(\,\cdot\,)=(T(\,\cdot\,)\xi,\xi)$. The
$\psi$'s are called the {\it spherical functions.\/} These are certain
$K$--biinvariant functions on $G$, which can be converted (via
restriction to the subgroup $U(\infty)\times\{e\}\subset G$) to
certain central functions $\chi$ on $U(\infty)$. The functions $\chi$
thus obtained are called the {\it characters\/} of the group
$U(\infty)$. The correspondence $T\leftrightarrow\chi$ makes it
possible to employ both languages, that of spherical
representations and that of characters, each of which has its own
merits. To {\it irreducible\/} representations $T$ correspond {\it
extreme\/} characters $\chi$ (i.e., extreme points in the convex set
of all characters). Irreducible spherical representations of
$(G,K)$ and extreme characters of $U(\infty)$ admit a complete
description. They depend on countably many continuous parameters.  

\subhead (d) How to get ``natural'' representations \endsubhead
Set $G(N)=U(N)\times U(N)$ and let $K(N)$ be the diagonal in $G(N)$,
$N=1,2,\dots$\,. The homogeneous space $G(N)/K(N)$ can be identified
with the group space $U(N)$, and then the action of $G(N)$ on
$G(N)/K(N)$ turns into the two--sided action of $U(N)$ on itself. Let
$\Reg_N$ denote the quasi--regular representation of $G(N)$ in
$L^2(G(N)/K(N))$. This is nothing else than the biregular
representation of the compact group $U(N)$ whose decomposition is
determined by the Peter--Weyl theorem.  

We seek for a counterpart $\Reg_\infty$ of the representations
$\Reg_N$ for the pair $(G,K)$. As was mentioned above, there is no good
measure on $G/K$, hence no space $L^2(G/K)$. To overcome this
difficulty, two general recipes are known. First, to embed
representations $\Reg_N$ in each other and then to define
$\Reg_\infty$ as the inductive limit representation
$\varinjlim_{N\to\infty}\Reg_N$. Second, to embed $G/K$ into an 
appropriate $G$--space $\overline{G/K}$ possessing an invariant (or
quasiinvariant) measure $m$ and then to realize $\Reg_\infty$ in
$L^2(\overline{G/K},m)$.  

A realization of any of these recipes is by no means an automatic
exercise: 
 
A subtle point of the inductive limit construction is that there are
many different embeddings $\iota_N:\Reg_N\to\Reg_{N+1}$ for each $N$,
and the limit representation $\Reg_\infty$ highly depends of the chain
$\{\iota_N\}$ chosen. Trying all the possible chains $\{\iota_N\}$
one gets too many limit representations, so that a fine selection
rule must be imposed. 

As for the second way, we have to guess what the ambient space $\overline{G/K}$
should be. With $\overline{G/K}$ specified, we next have to find good measures $m$
(and moreover, to select some 1--cocycles that are also involved in the
construction).   

We employ both methods and finally get a family $\{T_{zw}\}$ of
representations depending on two complex parameters $z,w$ such that
$\Re(z+w)>-\tfrac12$. We believe that the representations $T_{zw}$ are
``natural'' objects of harmonic analysis. 

The space $\overline{G/K}$ is constructed as follows. We define certain
projections of the group spaces $U(N)\to U(N-1)$, where
$N=2,3,\dots$, and then take the projective limit space $\frak
U=\varprojlim U(N)$. This is our $\overline{G/K}$. 

\subhead (e) Gelfand--Tsetlin graph and coherent systems \endsubhead
The Gelfand--Tsetlin graph is a convenient tool for writing
characters of $U(\infty)$. The vertices of the graph symbolize the
irreducible representations of various groups $U(N)$ while the edges
code the inclusion relations between irreducible representations of
$U(N)$ and $U(N+1)$. By $\GT_N$ we denote the subset of
vertices corresponding to irreducibles of $U(N)$; elements of $\GT_N$
are identified with dominant weights for $U(N)$, i.e., these are
$N$--tuples $\la=(\la_1\ge\dots\ge\la_N)$ of integers.  

Given a character $\chi$ of $U(\infty)$, we can expand its restriction
to the subgroup $U(N)$ into a convex combination of the functions
$\chi^\la(\,\cdot\,)/\chi^\la(e)$, where $\chi^\la$ stands for the
irreducible character (in the conventional sense) of $U(N)$, indexed
by $\la\in\GT_N$. The coefficients $P_N(\la)$ of this expansion
determine a probability distribution $P_N$ on the discrete set
$\GT_N$. 

In this way, we get a bijection
$\chi\leftrightarrow\{P_N\}_{N=1,2,\dots}$ between characters $\chi$
and certain sequences $\{P_N\}$ of probability distributions. These
sequences are called {\it coherent systems\/}, because, for any $N$,
the distributions $P_N$ and $P_{N+1}$ are connected by a certain
``coherency relation''. 

\subhead (f) The characters $\chi_{zw}$ and their analytic
continuation \endsubhead
Having constructed the representations $T_{zw}$, we proceed to the
corresponding characters $\chi_{zw}$. We evaluate them in terms of
the associated coherent systems $\{P_N\}$. The expressions that we
get for $P_N$ make sense, via analytic continuation, for a larger set
of parameters. 

Specifically, let $\zw$ be the coordinates in $\C^4$, and let $\Cal
D$ be the open halfspace in $\C^4$ determined by the inequality
$\Re(z+w+z'+w')>-1$. We exhibit an ``admissible subset''
$\Dadm\subset\Cal D$ of real dimension 4, containing as a proper
subset all the quadruples $(\zw)\in\Cal D$ with $z'=\bar z$, 
$w'=\bar w$, and such that for any quadruple $(\zw)\in\Dadm$, the
following formulas provide a coherent system:
$$
\gather
P_N(\la\mid \zw)
=(S_N(\zw))^{-1}
\times\, \prod_{1\le i<j\le N}\frac{(\la_i-\la_j+j-i)^2}{(j-i)^2}\\
\times\,\prod_{i=1}^N
\frac1{\Gamma(z-\la_i+i)\Gamma(z'-\la_i+i)
\Gamma(w+N+1+\la_i-i)\Gamma(w'+N+1+\la_i-i)}\,,
\endgather
$$
where $N=1,2,\dots$, $\la$ ranges over $\GT_N$, and $S_N(\zw)$ is a
normalization constant:
$$
S_N(\zw)=\prod_{i=1}^N
\frac{\Gamma(z+z'+w+w'+i)}
{\Gamma(z+w+i)\Gamma(z+w'+i)\Gamma(z'+w+i)\Gamma(z'+w'+i)\Gamma(i)}\,.
$$

Thus, to any $(\zw)\in\Dadm$, we may assign a character $\chi_\zw$.
The initial characters $\chi_{zw}=\chi_{z,w,\bar z,\bar w}$ form the
``principal series'', and the remaining characters belong to its
analytic continuation. The whole picture resembles the conventional
principal, complementary and degenerate series for semi--simple Lie
groups, so that we even employ this terminology. However, in our
context, these concepts refer not to (generically) irreducible
representations, as in the conventional context, but to highly
reducible ones.  

Note that the construction of the coherent systems implies a curious
summation formula, 
$$
\sum\Sb \la_1\ge\dots\ge\la_n\\ 
\la_1,\dots,\la_N\in\Bbb Z\endSb 
P_N(\la\mid\zw)=1.
$$
In the simplest case $N=1$ it looks as 
$$
\gather
\sum_{k\in\Bbb Z}
\frac1{\Gamma(z-k+1)\Gamma(z'-k+1)
\Gamma(w+k+1)\Gamma(w'+k+1)}\\
=\frac{\Gamma(z+z'+w+w'+1)}
{\Gamma(z+w+1)\Gamma(z+w'+1)\Gamma(z'+w+1)\Gamma(z'+w'+1)}\,,
\endgather
$$
which is equivalent to a classical identity due to Dougall, see
\cite{AAR, Chapter 2, Theorem 2.8.2 and Exercise 42(b)}, \cite{Er, \S1.4}. 

\subhead (g) Abstract theorems on spectral measures \endsubhead
Let us return now to extreme characters. There is a bijective
correspondence $\chi^{(\om)}\leftrightarrow\om$ between extreme
characters and points $\om$ of an infinite--dimensional ``region''
$$
\Om\,\subset\,\R^{4\infty+2}=\R^\infty\times\R^\infty\times
\R^\infty\times\R^\infty\times\R\times\R.
$$
For the functions $\chi^{(\om)}(U)$, where $U\in U(\infty)$, there is a
remarkable explicit formula due to Voiculescu (see \tht{1.2} below). 

Our first ``abstract'' theorem says that for any character $\chi$,
here exists a unique probability measure $P$ on $\Om$ such that
$$
\chi(U)=\int_\Om \chi^{(\om)}(U)\,P(d\om), \qquad U\in U(\infty).
$$
We call $P$ the {\it spectral measure\/} of $\chi$. Conversely, any
probability measure $P$ on $\Om$ is a spectral measure for a certain
$\chi$. This result is a refinement of a theorem due to Voiculescu
\cite{Vo}. It provides a nice general description of the whole
set of characters. 

The next result shows that spectral measures can be, in
principle, computed. Specifically, we define, for any $N=1,2,\dots$,
an embedding $\GT_N\hookrightarrow\Om$ such that the image of $\GT_N$
looks as a discrete approximation of $\Om$, which becomes more and
more exact as $N\to\infty$.  

Now let $\chi$ be an arbitrary character, $\{P_N\}$ be the
corresponding coherent system, and $P$ be the spectral measure of
$\chi$. Then, according to the second ``abstract'' theorem, we have
$P=\lim_{N\to\infty}P_N$, where we identify $P_N$ with its
pushforward under the embedding $\GT_N\hookrightarrow\Om$.  

\subhead (h) The problem of harmonic analysis on $U(\infty)$
\endsubhead 
Now we are in a position to state this problem explicitly: 

Let $\chi=\chi_\zw$, where $(\zw)\in\Dadm$, and let $\{P_N\}$ be the
corresponding coherent system. Recall that
$P_N=P_N(\,\cdot\,\mid\zw)$ is a probability measure on $\GT_N$, 
given by the explicit formula above. Let us carry over $P_N$ to the
space $\Om$. Then the problem consists in evaluating the limit of the
measures $P_N$ in the ambient space $\Om$ as $N\to\infty$. 

(By virtue of the ``abstract'' theorems above, the limit always exists
and coincides with the spectral measure.)

The solution to the problem is presented in the next paper,
\cite{BO4}. 

\subhead (i) Connections with infinite random matrices \cite{BO3}
\endsubhead 
Let $\HH$ be the space of all infinite Hermitian matrices and let
$U(\infty)$ act on $\HH$ by conjugations. There is a parallelism:

\medskip    
{\settabs 2 \columns
\+ Characters $\chi$ of $U(\infty)$. &  
Invariant probability measures $M$ on $\HH$.  \cr
\+Extreme characters $\chi^{(\om)}$,&
Ergodic invariant measures $M^{(\om)}$, \cr
\+indexed by points $\om$ &
indexed by points $\om$ \cr
\+of a region $\Om\subset\R^{4\infty+2}$. & 
of a region $\Om\subset\R^{2\infty+2}$. \cr
\+Decomposition $\chi=\int_\Om\chi^{(\om)} P(d\om)$, &
Decomposition $M=\int_\Om M^{(\om)} P(d\om),$ \cr
\+where $P$ is a probability measure &
where $P$ is a probability measure \cr
\+on $\Om\subset\R^{4\infty+2}$. &
on $\Om\subset\R^{2\infty+2}$. \cr
\+Distinguished characters $\chi=\chi_{zw}$, &
Distinguished invariant measures $M=m^{(s)}$, \cr
\+where $z,w\in\C$, $\Re(z+w)>-\tfrac12$. &
where $s\in\C$, $\Re s>-\tfrac12$.\cr 
\+The problem of computing &
The problem of computing \cr
\+the spectral measures $P$ &
the spectral measures $P$ \cr
\+for distinguished characters $\chi$. &
for distinguished measures $M$. \cr}
\medskip
\noindent And so on.

In the paper \cite{BO3}, we constructed the distinguished measures
$m^{(s)}$ and computed their spectral measures. The whole theory of
$U(\infty)$--invariant measures can be viewed as a simplified version
of the theory of characters. 

On the other hand, the results of \cite{BO3} are directly used
in the present paper: the space $\frak U$ mentioned above in
subsection (d) can be identified (within a negligible subset) with
the space $\HH$, and the measures $m$ involved in the construction of
the representations $T_{zw}$ are nothing else than the measures $m^{(s)}$. 

\subhead (j) Characters of $S(\infty)$ \endsubhead
As was mentioned above, the group $S(\infty)$, the union of finite
symmetric groups $S(n)$, is another fundamental example of an
inductive limit group. The representation theories of the both
groups, $U(\infty)$ and $S(\infty)$, reveal deep analogies and links.
Of course, in the classical representation 
theory, it is well known that representations of the groups $S(n)$
and $U(N)$ are connected by the Schur--Weyl duality. But in the
case of infinite dimension the connections between the symmetric and
unitary groups turn out to be much closer.

The problem of harmonic analysis on $S(\infty)$ was stated in
\cite{KOV} and then further developed in a cycle of papers
\cite{P.I-V}, \cite{Bor1-2}, \cite{BO1-2}. The construction
\cite{KOV} of the ``generalized regular 
representations'' $T_z$ served as a guiding example for the
construction of the representations $T_{zw}$ in the present paper.
Further, the experience of our work on the spectral measures for
$S(\infty)$ helped us very much in the work \cite{BO4}.

\subhead (k) Generalization to other $(G,K)$--pairs. Works of
Pickrell \cite{Pi} and Neretin \cite{Ner3} \endsubhead 
The $(G,K)$--pair $(U(\infty)\times U(\infty), \operatorname{diagonal}
U(\infty))$ is a representative of the family of ten $(G,K)$--pairs,
which come from the ten classical series of compact Riemannian
symmetric spaces. This family is a natural framework for developing
representation theory (see \cite{Ol1}, \cite{Ol3}, \cite{Ner1}) and,
in particular, harmonic analysis. In the pioneer work \cite{Pi},
Pickrell considered 
the pair $G=\varinjlim U(2N)$, $K=U(N)\times U(N)$, which corresponds to
the series of complex Grassmanians $U(2N)/U(N)\times U(N)$. He
constructed the ambient space $\overline{G/K}$ as a projective limit of
Grassmanians and defined a family of measures $m$ which give rise
to ``natural'' representations. Pickrell's construction was the
starting point of \cite{KOV} and of the present paper. In the recent
paper by Neretin, \cite{Ner3}, Pickrell's construction is carried over
to all ten pairs. 

Note that, as compared with Pickrell's results, our construction of the
representations $T_{zw}$ incorporates a few new observations. 

First, following \cite{KOV}, we introduce a complex parameter instead
of a real one (the idea is to employ a wider family of cocycles for
an action of the group $G$). It is worth noting that the same
generalization makes sense for all ten  pairs $(G,K)$ mentioned
above. 

Second, for our pair $(G,K)$, it is actually
possible to introduce {\it two\/} complex parameters $z,w$. This
(nonevident) fact was prompted by Neretin's results on Hua--type 
integrals, see \cite{Ner2}, \cite{Ner3}. 

Third, we observe that the construction of the representations
$T_{zw}$ makes sense for all $z,w\in\C$. But for $\Re z+\Re w\le-1/2$
we get representations without a distinguished $K$--invariant vector.
It would be interesting to study these representations. 

\subhead (l) Acknowledgment \endsubhead
This paper is part of a joint project with Alexei Borodin. I am very
grateful to him and also to Sergei Kerov, Yurii Neretin, and Anatoly
Vershik for numerous discussions.  

\head 1. Characters  \endhead

\example{Definition 1.1} 
Let $K$ be a topological group. By a {\it character\/} of $K$ we 
mean any continuous complex--valued function $\chi$ on $K$ satisfying
the following three conditions: 
 
(i) $\chi$ is central, i.e., constant on conjugacy classes; 
 
(ii) $\chi$ is positive definite, i.e., for any finite collection
$g_1,\dots,g_n$ of elements of $K$, the $n\times n$ matrix
$[\chi(g_j^{-1}g_i)]_{1\le i,j\le n}$ is Hermitian and nonnegative; 
 
(iii) $\chi$ is normalized at the unit element, i.e.,
$\chi(e)=1$. 

Let $\Cal X(K)$ denote the set of the characters of $K$. Evidently,
$\Cal X(K)$ is a convex set. Its extreme points are called {\it extreme\/}
(or {\it indecomposable\/}) characters. \footnote{This terminology
differs from that used in \cite{Th1}, \cite{Th2}, \cite{VK1},
\cite{VK2}, \cite{Vo}. In those papers, only extreme characters were
considered, and they were called simply characters.}
\endexample

\example{Example 1.2} Let $K$ be a finite group or, more generally, a
compact separable group, and let $\widehat K$ be its dual space,
i.e., the set of equivalence classes of 
irreducible representations. Then $\widehat K$ is a finite or a
countably infinite set.  Let $\la$ range over $\widehat K$ and $\chi^\la$
denote the irreducible 
character corresponding to $\la$ (i.e., the trace of an irreducible
representation in the class $\la$). The extreme characters of $K$ in the
sense of Definition 1.1 are exactly the {\it normalized irreducible
characters}
$$
\wt\chi^\la(g)=\frac{\chi^\la}{\chi^\la(e)}\,, 
\qquad \la\in\widehat K, \tag1.1
$$
while general characters $\chi\in \Cal X(K)$ are convex combinations of
the extreme ones,
$$
\chi=\sum_{\la\in\widehat K}P(\la)\wt\chi^\la, \qquad 
P(\la)\ge0, \quad \sum_{\la\in\widehat K}P(\la)=1.
$$
Thus, the set $\Cal X(K)$ is a simplex with vertices indexed by
$\la\in\widehat K$.  
\endexample

In the example above the 
set $\Cal X(K)$ is large enough to separate conjugacy classes. If $K$ is
not compact it may well happen that $\Cal X(K)$ is very small (e.g., is
exhausted by the function $\chi\equiv1$). Actually, the class of
groups $K$ with $\Cal X(K)$ large enough is rather restricted. However, it
includes important examples leading to a rich theory. 

In the present paper we take as $K$ the {\it infinite--dimensional
unitary group,\/} which is defined as  
$$
U(\infty)=\bigcup_{N\ge1} U(N),
$$
where  $U(N)$ is
the group of $N\times N$ unitary matrices. The embedding 
$U(N)\hookrightarrow U(N+1)$ is defined as follows: we identify $U(N)$
with the subgroup 
in $U(N+1)$ fixing the $(N+1)$st basis vector. Equivalently,
$U(\infty)$ is the group of infinite unitary matrices 
$U=[U_{ij}]$, $i,j=1,2,\dots$, with finitely
many matrix entries $U_{ij}$ distinct from $\delta_{ij}$. 

We equip $U(\infty)$ with the inductive limit topology (i.e., a
function on $U(\infty)$ is continuous if its restriction to each
subgroup $U(N)$ is continuous). Note that $U(\infty)$ is {\it not\/}
locally compact. 

First, let us describe the conjugacy classes in $U(\infty)$.  
Recall that the conjugacy classes in $U(N)$ are parameterized by spectra of
unitary matrices, that is, by 
unordered $N$-tuples $u_1,\dots,u_N$ of complex numbers with modulus
1. In this notation, the embedding $U(N)\hookrightarrow U(N+1)$ is
described as $(u_1,\dots,u_n)\mapsto(u_1,\dots,u_N,1)$. Thus, the
conjugacy classes in 
$U(\infty)$ can be parameterized by countable collections of complex
numbers $(u_1,u_2,\dots)$ such that $|u_i|=1$ and only finitely many of
$u_i$'s are different from 1; the ordering of $u_i$'s 
is unessential. As a representative of a conjugacy class indexed by
$(u_1,u_2,\dots)$ one can take the diagonal matrix
$\operatorname{diag}(u_1,u_2,\dots)$. 

To describe the extreme characters of $U(\infty)$ we need some
notation. 

Let $\R^\infty$ denote the product of countably many
copies of $\R$, and set
$$
\R^{4\infty+2}=\R^\infty\times\R^\infty\times\R^\infty\times\R^\infty
\times\R\times\R.
$$
Let $\Om\subset\R^{4\infty+2}$ be the subset of sextuples
$$
\om=(\al^+,\be^+;\al^-,\be^-;\de^+,\de^-)
$$ 
such that 
$$
\gather
\al^\pm=(\al_1^\pm\ge\al_2^\pm\ge\dots\ge 0)\in\R^\infty,\quad
\be^\pm=(\be_1^\pm\ge\be_2^\pm\ge\dots\ge 0)\in\R^\infty,\\
\sum_{i=1}^\infty(\al_i^\pm+\be_i^\pm)\le\de^\pm, \quad
\be_1^++\be_1^-\le 1.
\endgather
$$
Set
$$
\ga^\pm=\de^\pm-\sum_{i=1}^\infty(\al_i^\pm+\be_i^\pm)
$$
and note that $\ga^+,\ga^-$ are nonnegative.

To any $\om\in\Om$ we assign a function $\chi^{(\om)}$ on $U(\infty)$:
$$
\chi^{(\om)}(U)=
\prod_{u\in\operatorname{Spectrum(U)}}
\left\{e^{\ga^+(u-1)+\ga^-(u^{-1}-1)}
\prod_{i=1}^\infty\frac{1+\be_i^+(u-1)}{1-\al_i^+(u-1)}
\,\frac{1+\be_i^-(u^{-1}-1)}{1-\al_i^-(u^{-1}-1)}\right\}\,.
\tag1.2
$$
Here $U$ is a matrix from $U(\infty)$ and $u$ ranges over the set of
its eigenvalues. All but finitely many $u$'s equal 1, so that the
product over $u$ is actually finite. The product over $i$ is
convergent, because the sum of the parameters is finite. 

\proclaim{Theorem 1.3} The functions $\chi^{(\om)}$, where $\om$ ranges
over $\Om$, are exactly the extreme characters of the group
$U(\infty)$. 
\endproclaim

The coordinates $\al^\pm_i$, $\be^\pm_i$, and $\ga^\pm$ (or $\de^\pm$)
are called the {\it Voiculescu parameters\/} of the extreme character
$\chi^{(\om)}$. Theorem 1.3 is similar to
Thoma's theorem which describes the extreme characters of the
infinite symmetric group, see \cite{Th1}, \cite{VK1}, \cite{Wa}, \cite{KOO}.

The study of extreme characters of the group $U(\infty)$ was
initiated by Voiculescu \cite{Vo}. He proved (among other things)
that all functions \tht{1.2} are actually extreme characters. 
The fact that Voiculescu's list is exhaustive was established later by 
Boyer \cite{Boy1} and by Vershik--Kerov \cite{VK2}. They independently 
pointed out that Theorem 1.3 is implied by old Edrei's result
\cite{Ed} about two--sided totally positive sequences. 

On the other hand, Vershik and Kerov outlined in \cite{VK2} another
approach to Theorem 1.3, based on the idea of approximating 
extreme characters by normalized irreducible characters of the
compact groups $U(N)$. The same idea was employed in  
their previous work \cite{VK1} on the infinite symmetric group
$S(\infty)$. 

A detailed proof of Theorem 1.3 by Vershik--Kerov's asymptotic method
was given much later by Okounkov--Olshanski \cite{OkOl} (this paper also
contains more general results). In a special case, a detailed proof
was earlier given by Boyer \cite{Boy2}.  

Here are some comments to Voiculescu's formula \tht{1.2}.

\example{Remark 1.4} The simplest extreme characters are of the form
$\chi(U)=\det^k(U)$, where $k\in\Z$. The
corresponding parameters are as follows: all of them are
equal to zero except the first $|k|$ coordinates in $\beta^+$ (if $k>0$)
or in $\be^-$ (if $k<0$), which are equal to 1. 
\endexample

\example{Remark 1.5} Given a character $\chi$, define 
the function $\chi\otimes\det^k(\,\cdot\,)$ as the pointwise product
$\chi(U)\det^k(U)$). Then $\chi\otimes\det^k(\,\cdot\,)$ is a
character, too. If $\chi$ is extreme then
$\chi\otimes\det^k(\,\cdot\,)$ is extreme. In terms of Voiculescu's
parameters, tensoring with $\det(\,\cdot\,)$ reduces to 
$$
\be^+\mapsto(1-\be^-_1,\be^+_1,\be^+_2,\dots), \qquad
\be^-\mapsto(\be^-_2,\be^-_3,\dots).
$$
\endexample

\example{Remark 1.6} Here is a comment to the condition
$\be^+_1+\be^-_1\le1$. Even if this condition is dropped,
Voiculescu's formula still defines an extreme character. However,
using the identity
$$
[1+b^+(u-1)][1+b^-(u^{-1}-1)]=
[1+(1-b^-)(u-1)][1+(1-b^+)(u^{-1}-1)],
$$
we can always modify the beta parameters so that the condition
$\be^+_1+\be^-_1\le1$ will be satisfied. With this condition imposed,
no freedom to change the parameters remains: if $\om_1\ne\om_2$ then
$\chi^{(\om_1)}\ne\chi^{(\om_2)}$, see \cite{OkOl, \S5, Step 3}.  
\endexample

\example{Remark 1.7} Let $SU(\infty)$ denote the subgroup of the
matrices $U\in U(\infty)$ with determinant 1. I.e., $SU(\infty)$ is
the union of the groups $SU(N)$. Restricting any character to
$SU(\infty)$ one gets a character of this group, and extreme
characters remain extreme. However, it may well happen that
$\chi^{(\om_1)}\mid_{SU(\infty)}$ equals
$\chi^{(\om_2)}\mid_{SU(\infty)}$ for certain $\om_1\ne\om_2$. As 
parameters for the characters $\chi^{(\om)}\mid_{SU(\infty)}$ one may take
$\al^\pm$, $\ga^\pm$, and a multiset $B$ in
$[-\tfrac12,\tfrac12]$, which is obtained by mixing together $\be^+$
and $\be^-$. Specifically,
$$
B=\{\be^+_i-\tfrac12\}_{i=1,2,\dots}\cup
\{-\be^-_j+\tfrac12\}_{j=1,2,\dots}.
$$
Given a point $b\in B$, it is no longer possible to decide whether it
comes from a coordinate of $\be^+$ or of $\be^-$. 
\endexample

\example{Remark 1.8} On the set of all characters, there is a natural
operation: pointwise conjugation. For extreme characters, it reduces
to the transposition
$(\al^+,\be^+,\de^+)\leftrightarrow(\al^-,\be^-,\de^-)$.  
\endexample

\example{Remark 1.9} One can check that the functions \tht{1.2}
separate conjugacy classes in $U(\infty)$.
\endexample

\example{Remark 1.10} The reader might notice that each extreme
character $\chi^{(\om)}$ is a multiplicative function with respect to a
natural product on the set of conjugacy classes of $U(\infty)$: taking
disjoint union of two collections of eigenvalues. Such a
multiplicativity property holds for many other groups and can be
established independently of the classification of extreme
characters, see the survey \cite{Ol4}.  
\endexample

\head 2. Spherical and admissible representations \endhead

There are two ways to establish a connection between characters and
unitary representations:

First, extreme characters of a group $K$ parameterize its finite
factor representations, i.e., unitary representations 
generating finite von Neumann factors. See, e.g., \cite{Th2}.

Second, extreme characters of $K$ parameterize irreducible
spherical representations of the pair $(G,\diag K)$, where
$G=K\times K$ and $\diag K$ is the diagonal subgroup in $G$. See
\cite{Ol1} and \cite{Ol3, Theorem 24.5}. In these papers, it is
also explained how to establish the relationship between both kinds
of representations directly.  

In the present paper we are dealing with spherical representations.
Below we review basic facts concerning the correspondence
between characters and spherical representations, specialized
to the particular case $K=U(\infty)$. 
For proofs of the claims stated below we refer to \cite{Ol3}.

$\bullet$ Set $G=U(\infty)\times U(\infty)$ and $K=\diag U(\infty)$.
Let $\Psi$ denote the set of continuous functions $\psi$ on the group
$G$ that are $K$-biinvariant, positive definite, and normalized at
the unit element. Such functions will be called {\it spherical
functions.\/} The set $\Psi$ is convex.

$\bullet$ There exists a natural bijective correspondence
$\chi\leftrightarrow \psi$ between characters $\chi\in \Cal X(U(\infty))$
and spherical functions $\psi\in\Psi$, which is defined as follows
$$
\chi(U)=\psi(U,1), \qquad 
\psi(U_1,U_2)=\chi(U_1U_2^{-1}), \qquad
U,U_1,U_2\in U(\infty).
$$

$\bullet$ The correspondence $\chi\leftrightarrow\psi$ is an
isomorphism of convex sets, so that extreme characters exactly
correspond to extreme spherical functions. 

$\bullet$ Any extreme function $\psi\in\Psi$ is also extreme in a wider
convex set, which is formed by all (not necessarily $K$-biinvariant)
positive definite normalized functions on $G$. 

$\bullet$ By a {\it spherical representation\/} of $(G,K)$ we mean
any pair $(T,\xi)$, where $T$ is a unitary representation of $G$ and
$\xi$ is a fixed cyclic unit $K$-invariant vector in the Hilbert
space of $T$. The vector $\xi$ is called the {\it spherical vector.}

$\bullet$ There is a bijective correspondence 
$\psi\leftrightarrow (T,\xi)$ between spherical functions and 
(equivalence classes of) spherical representations, defined by 
$$
\psi(g)=(T(g)\xi,\xi), \qquad g=(U_1,U_2)\in G.
$$
Equivalently, there is a bijection $\chi\leftrightarrow(T,\xi)$
between characters and (equivalence classes of) spherical
representations, defined by
$$
\chi(U)=(T(U,1)\xi,\xi), \qquad U\in U(\infty).
$$

$\bullet$ Under this bijection, extreme spherical functions (or
extreme characters) exactly correspond to irreducible spherical
representations. 

$\bullet$ In an irreducible spherical representation, the subspace of
$K$-invariant vectors has dimension 1. 

\medskip

Combining these claims with Theorem 1.3 we get the following

\proclaim{Corollary 2.1} Irreducible spherical representations of the
pair $(G,K)$, where $G=U(\infty)\times U(\infty)$ and $K=\diag
U(\infty)$, are parameterized by the points $\om\in\Om$, where the
region $\Om\subset\R^{4\infty+2}$ is defined in \S1.
\endproclaim

An explicit realization of the irreducible spherical
representations is given in \cite{Ol3}, \cite{Ol2}.

The class of spherical representations is part of a wider class of
unitary representations, which is defined as follows.

\example{Definition 2.2} Let $G,K$ be as above and let $K_n$ denote
the subgroup of $K\simeq U(\infty)$ constituted by matrices of the form 
$\bmatrix 1_n & 0\\ 0 & *\endbmatrix$. A unitary representation $T$ of the
group $G$ is called an {\it admissible representation\/} of $(G,K)$ if
the subspace of $K_n$-invariant vectors (with varying $n$) is dense
in the Hilbert space of $T$. 
\endexample 

It is readily proved that any spherical representation is admissible.
A (conjecturally, exhaustive) list of irreducible admissible
representations, together with their explicit realization, can be
found in \cite{Ol3}, \cite{Ol2}. Each representation from this list
is specified by a point $\om\in\Om$ and some additional
discrete data.  
It is worth noting that, although $G$ is the product of two copies of
$U(\infty)$, generic irreducible admissible representations of
$(G,K)$ are not tensor products of two irreducible representations 
of $U(\infty)$.

\head 3. The representations $T_{z,w}$ \endhead

In this section we construct a family of representations of the group
$G=U(\infty)\times U(\infty)$ depending on 2 complex parameters $z,w$.

Let us abbreviate
$$
G(N)=U(N)\times U(N)\subset G, \qquad K(N)=\diag U(N)\subset K,
\quad N=1,2,\dots\,.
$$
We consider $U(N)$ as the homogeneous space $K(N)\setminus G(N)$ with
the following right action of $G(N)$:
$$
(U,\,(U_1,U_2))\mapsto U_2^{-1}UU_1\,, \qquad U\in U(N), 
\quad (U_1,U_2)\in G(N). 
$$

Let $N\ge2$. Given a unitary matrix $U\in U(N)$,
write it in the block form corresponding to the partition
$N=(N-1)+1$: 
$$
U=\bmatrix A & B\\ C & D\endbmatrix,
$$
so that $D=U_{NN}$ is a complex number while $A$ is a $(N-1)\times(N-1)$
matrix, and set
$$
p_N(U)=\cases A-B(1+D)^{-1}C, & D\ne-1,\\
A, & D=-1. \endcases
$$

\proclaim{Lemma 3.1} The map $p_N$ defines a projection $U(N)\to
U(N-1)$, which commutes with the action of $G(N-1)$, that is, with
left and right shifts by elements of $U(N-1)$.
\endproclaim 

We will call $p_N$ the {\it canonical projection.} The claim of the
lemma is by no means new (see, e.g., \cite{Ner3}). For
the reader's convenience we give a detailed proof below.

\demo{Proof} The fact that $p_N$ commutes with left and right shifts
by elements of $U(N-1)$ is evident. Next, if $D=-1$ then $B=0$,
$C=0$, and $A$ is a unitary matrix; therefore, in this case
$p_N(U)\in U(N-1)$. It remains to prove that if $D\ne-1$ then the
matrix $V=A-B(1+D)^{-1}C$ is unitary. 

Let $\Ga\subset\C^N\oplus\C^N$ be the graph of 
$\wt U=\bmatrix \phantom{-}A & \phantom{-}B\\ -C & -D\endbmatrix$, that is
$$
\Ga=\{x\oplus y\mid x\in\C^N, \, y\in\C^N, \, y=\wt Ux\}.
$$
Write
$$
x=x_1\oplus x_2, \qquad y=y_1\oplus y_2, \qquad
x_1, \,y_1 \in \C^{N-1}\,, \quad x_2,\, y_2 \in\C^1\,.
$$
Intersect $\Ga$ with the hyperplane $x_2=y_2$ and denote by $\Ga_1$
the image of this intersection under the projection
$$
\C^N\oplus\C^N\to\C^{N-1}\oplus\C^{N-1}, \qquad
x\oplus y\mapsto x_1\oplus y_1\,.
$$
Clearly, $\Ga_1$ is described by the linear equations
$$
\bmatrix y_1\\y_2\endbmatrix=
\bmatrix \phantom{-}A & \phantom{-}B\\ -C & -D\endbmatrix\,
\bmatrix x_1\\x_2\endbmatrix\,, \qquad 
y_2=x_2.
$$
Or, equivalently,
$$
y_1=Ax_1+Bx_2\,, \qquad
y_2=-Cx_1-Dx_2\,, \qquad
y_2=x_2\,.
$$
Excluding $x_2$ and $y_2$ from this system we conclude that $\Ga_1$
coincides with the graph of $V$. 

Since $\wt U$ is unitary, we have $(x,x)=(y,y)$ for any 
$x\oplus y\in\Ga$, or, equivalently,
$(x_1,x_1)+|x_2|^2=(y_1,y_1)+|y_2|^2$. Therefore, $x_2=y_2 $ implies
$(x_1,x_1)=(y_1,y_1)$, which means that $V$ is unitary. 
\qed
\enddemo

\example{Remark 3.2} The above interpretation of the projection $p_N$
in terms of graphs of operators also works in the exceptional case $D=-1$.
More generally, it can be used to describe the projection 
$p_{M,N}=p_{M+1}\circ\dots\circ p_N$ from $U(N)$ to $U(M)$ for any
$M<N$. To do this we have to write the matrices $U\in U(N)$ in the
block form corresponding to the partition $N=M+(N-M)$ and then
repeat the same argument.
\endexample

\example{Remark 3.3} The projection $p_N$ is continuous on the open
subset of $U(N)$ consisting of matrices $U$ with $U_{NN}\ne-1$ but
not on the whole group $U(N)$. Indeed, the first claim is evident. To
demonstrate discontinuity, take $N=2$ and remark that 
$$
p_2\left(\bmatrix \cos\varphi & \sin\varphi \\
-\sin\varphi & \cos\varphi\endbmatrix\right)=
\cases 1, & \cos\varphi\ne-1\\ -1, & \cos\varphi=-1.\endcases
$$
\endexample

\example{Remark 3.4} The projection $p_{M,N}:U(N)\to U(M)$ mentioned
in Remark 3.2 is connected with characteristic functions in the sense
of Liv\u{s}ic et al. Write $U\in U(N)$ in the block form
corresponding to the
partition $N=M+(N-M)$: 
$$
U=\bmatrix A & B\\ C & D\endbmatrix,
$$
where $A$ is of the size $M\times M$. The matrix--valued function
$$
\phi_U(\zeta)=A+\zeta B(1-\zeta D)^{-1}C
$$ is called the {\it characteristic function\/} of $U$: it is an
important invariant of the action of the subgroup 
$$
\left\{g\in U(N)\mid g=\bmatrix 1_M & 0\\ 0 & *\endbmatrix\right\} 
\simeq U(N-M)
$$
on $U(N)$ by conjugations. The function $\phi_U(\zeta)$ is an inner function:
$\Vert\phi_U(\zeta)\Vert<1$ in the open disk $|\zeta|<1$ while
$\phi_U(\zeta)\in U(M)$ on the
boundary $|\zeta|=1$.  See \cite{Ner1, Appendix E} for further
information and references to original papers. We have 
$p_{M,N}(U)=\phi_U(-1)$. 
\endexample

Although $p_N$ is not continuous on the whole $U(N)$, it is clearly a
Borel map. Hence, given a Borel measure on $U(N)$, we may project it
to $U(N-1)$ by means of $p_N$. 

In particular, let $\mu_N$ denote the normalized Haar measure on
$U(N)$. Then Lemma 3.1 implies that $p_N$ takes $\mu_N$ to an
invariant probability measure on $U(N-1)$, that is, to $\mu_{N-1}$.

Let $\DD=\{D\in\C\mid |D|\le1\}$ denote the closed
unit disk. Introduce a projection 
$$
\ep_N: U(N)\to\DD,\qquad
U\mapsto D=U_{NN}\,.
$$
This map is constant on double cosets modulo $U(N-1)$. Moreover, one
can prove that $\ep_N$ separates these cosets, i.e., if two matrices
have the same $(N,N)$--entry then they belong to the same double coset
modulo $U(N-1)$. 

Denote by $\nu_N$ the image of $\mu_N$ under $\ep_N$; this is a
probability measure on the disk $\DD$. 

\proclaim{Lemma 3.5} 
$$
\nu_N(dD)=\const\, (1-|D|^2)^{N-2}\Leb(dD),
$$
where $\Leb$ is the Lebesgue measure on $\DD$.
\endproclaim 

\demo{Proof} This claim is well known. It can be deduced from the
fact that $\nu_N$ coincides with the image of the uniform measure on
the $(2N-1)$--dimensional sphere
$\{\zeta\in\C^N\mid\,|\zeta_1|^2+\dots|\zeta_N|^2=1\}$ under the
projection $\zeta\mapsto\zeta_N\in\DD$. \qed
\enddemo

Following \cite{Ner3, \S1.2}, we combine $p_N$ and
$\ep_N$ into a single projection
$$
\wt p_N: U(N)\to U(N-1)\times\DD.
$$

\proclaim{Lemma 3.6} The image of $\mu_N$ under $\wt p_N$ is
$\mu_{N-1}\times\nu_N$.  
\endproclaim

\demo{Proof} Indeed, if $U\in U(N)$ and $\wt p_N(U)=(V,D)$ then, for
any $U_1,U_2\in U(N-1)$, we have $\wt
p_N(U_2^{-1}UU_1)=(U_2^{-1}VU_1,D)$. It follows that the
pushforward of $\mu_N$ under 
$\wt p_N$ splits into the direct product of $\mu_{N-1}$ with a
certain probability measure on $\DD$. Then the latter measure must
coincide with $\nu_N$. \qed
\enddemo

Set $H_N=L^2(U(N),\mu_N)$ and denote by $\Reg_N$ the
(bi)regular representation of $G(N)$ in $H_N$:
$$
(\Reg_N(g)f)(U)=f(U_2^{-1}UU_1), \qquad f\in H_N, \quad U\in U(N),
\quad g=(U_1,U_2)\in G(N).
$$
Next, introduce the subspace $H^0_N\subset H_N$ of functions
depending only on $\wt p_N(\,\cdot\,)$. 

\proclaim{Lemma 3.7} $H^0_N$ is a $G(N-1)$-invariant
subspace, and we have a natural isometry
$$
H^0_N\simeq H_{N-1}\otimes V_N, \qquad
V_N=L^2(\DD,\nu_N).
$$
\endproclaim

\demo{Proof} Indeed, this follows from the definitions and Lemma 3.6.
\qed
\enddemo
 
By Lemma 3.7, any unit vector $v\in V_N$ determines an embedding
$\Reg_{N-1}\to\Reg_N$, 
$$
H_{N-1}\ni f\mapsto f\otimes v\in H^0_N\subset H_N\,.
$$

Fix two complex numbers $z,w$. Our next aim is to assign a meaning to the
expression  
$$
f_{z,w\mid N}(U)=\det((1+U)^z(1+U^{-1})^w), \qquad U\in U(N).
$$
To do this we need a preparation which concludes in Definition 3.10. 

Let $\M$ denote the space of complex $N\times N$ matrices. For
$X\in\M$ let $\Re X=\frac12(X+X^*)$. We set
$$
\gather
\Mplus=\{X\in\M\mid \Re X>0\},\\
\MM=\{X\in\M\mid \Re X>-1\}=-1+\Mplus\,.
\endgather
$$
These are open domains in $\M$, isomorphic to a matrix wedge.

\proclaim{Lemma 3.8} For any fixed $z\in\C$, the expression
$(1+X)^z$ makes sense for $X\in\MM$ and defines a holomorphic map
$\MM\to GL(N,\C)$. Moreover, 
$$
(1+X)^z\cdot(1+X)^{z'}=(1+X)^{z+z'}, 
\qquad X\in\MM, \quad z,z'\in\C.
$$
\endproclaim

\demo{Proof} Setting $1+X=Y$ we must prove that $Y\mapsto Y^z$ is a
correctly defined holomorphic map $\Mplus\to GL(N,\C)$.

First, remark that the spectrum of any matrix $Y\in\Mplus$ lies in
the open half--plane $\C_+=\{\zeta\mid\Re\zeta>0\}$. Indeed, let
$\zeta$ be an 
eigenvalue of $Y$ and $v\ne0$ be any eigenvector corresponding to
$\zeta$. By the definition of $\Mplus$,
$$
0<((Y+Y^*)v,v)=\zeta+\bar\zeta,
$$
which means that $\zeta\in\C_+$. 

Now we can define $Y^z$ via the functional calculus:
$$
Y^z=\int_C \zeta^z(\zeta1-Y)^{-1}d\zeta,
$$
where $\zeta^z$ is correctly defined in $\C_+$ (we choose the branch
taking value 1 at $\zeta=1\in\C_+$) and $C$ is any simple contour in $\C_+$ 
containing the spectrum of $Y$. 

Clearly, $Y^z$ is nonsingular together with $Y$.

Now, to check the second claim it suffices to remark that it
obviously holds for matrices $X$ close to 0.
\qed 
\enddemo

\proclaim{Corollary 3.9} For any fixed $z\in\C$, the function
$X\mapsto\det((1+X)^z)$ is correctly defined and holomorphic in
$\MM$. Moreover,
$$
\det((1+X)^z)\det((1+X)^{z'})=\det((1+X)^{z+z'}), 
\qquad X\in\MM, \quad z,z'\in\C.
$$
\endproclaim 
\qed

Note that, denoting by $x_1,\dots,x_N$ the eigenvalues of $X$, we
have 
$$
\det((1+X)^z)=\prod_{k=1}^N(1+x_k)^z.
$$
The right--hand side makes sense, because $\Re x_k>-1$. On the
contrary, the expression $(\det(1+X))^z=(\prod(1+x_k))^z$ is
ambiguous. 

\example{Definition 3.10} Let $U(N)'\subset U(N)$ denote the set of
unitary matrices which do 
not have $-1$ as an eigenvalue. Note that $U(N)'=U(N)\cap\MM$. By
Corollary 3.9, the function 
$$
f_{z,w\mid N}(U)=\det((1+U)^z(1+U^{-1})^w)
$$ 
is well defined on $U(N)'$. 

Equivalently, denoting by $u_1,\dots,u_N$ the eigenvalues of $U$ and
assuming that none of them equals $-1$, we have
$$
f_{z,w\mid N}(U)=\prod_{k=1}^N(1+u_k)^z(1+\bar u_k)^w.
$$

On the complement of $U(N)'$, which is a
negligible set with respect to the Haar measure, we agree to continue
the function by 0. 
\qed 
\endexample

Thus, we have defined the function $f_{z,w\mid N}$
for any complex $z,w$ and each $N=1,2,\dots$ . 
In a few lemmas below we describe some special properties of these
functions which are used in the construction of the representations.

\proclaim{Lemma 3.11} Write arbitrary $N\times N$ matrices in block
form according to the partition $N=(N-1)+1$. The map 
$$
X=\bmatrix A&B\\ C&D\endbmatrix \to X_1=A-B(1+D)^{-1}C
$$
is correctly defined in the
domain $\MM$ and projects it onto the domain $\MMa$.

Consequently, $p_N$ maps $U(N)'$ onto $U(N-1)'$. 
\endproclaim

\demo{Proof} First of all, note that $X\in\MM$ implies $\Re D>-1$, so
that $1+D\ne0$ and the matrix $X_1$ is well defined. 

Next, we use the argument and notation of Lemma 3.1. Consider the graph
$\Ga$ of the operator 
$\bmatrix \phantom{-}A & \phantom{-}B\\ -C & -D\endbmatrix$. In terms
of $\Ga$, the condition $\Re X>-1$ means 
$$
(x_1,x_1)+|x_2|^2+\Re((x_1,y_1)-x_2\bar y_2)>0, 
\qquad x\oplus y\in\Ga.
$$
When $x_2=y_2$, this condition turns into 
$$
(x_1,x_1)+\Re(x_1,y_1)>0, \qquad x_1\oplus y_1\in\Ga_1,
$$
which means that $X_1\in\MMa$.

Together with Lemma 3.1 this implies that $p_N$ maps $U(N)'$ to
$U(N-1)'$. The surjectivity is evident. \qed
\enddemo

\proclaim{Lemma 3.12} For any $z\in\C$ and 
$X=\bmatrix A&B\\ C&D\endbmatrix\in\MM$ we have
$$
\det((1+X)^z)=\det((1+p_N(X))^z)\cdot(1+D)^z.
$$
\endproclaim

\demo{Proof} First of all, note that all the terms of this formula
make sense. Indeed, $(1+X)^z$ is correctly defined by Lemma 3.8,
$(1+p_N(X))^z$ 
is correctly defined by virtue of Lemma 3.8 and Lemma 3.11, and,
finally, $(1+D)^z$ is well defined, because $\Re D>-1$.

Next, we remark that these three terms are
holomorphic functions in $X\in\MM$. Hence, it suffices to prove the
formula when $X$ is near zero. Then we may interchange the symbol of
determinant and exponentiation, thus reducing the problem to the
identity 
$$
\det(1+X)=\det(1+p_N(X))\cdot(1+D),
$$
which is equivalent to  
$$
\det\bmatrix 1+A & B\\ C & 1+D\endbmatrix=
\det((1+A)-B(1+D)^{-1}C)\cdot(1+D).
$$
This is a special case of the well--known
identity for the determinant of a block matrix (of an arbitrary format),
$$
\det\bmatrix a & b\\ c & d\endbmatrix=
\det(a-bd^{-1}c)\,\det d.
$$
\qed
\enddemo

As a corollary we get the following formula: 
$$
f_{z,w\mid N}(U)=f_{z,w\mid N-1}(p_N(U))\cdot
(1+D)^z(1+\overline D)^w\,, \tag3.1
$$
where
$$
U\in U(N)', \qquad D=U_{NN}=\ep_N(U)\in\DD\setminus\{-1\}. 
$$

\proclaim{Lemma 3.13} The function $D\mapsto (1+D)^z(1+\overline D)^w$
is square 
integrable with respect to the measure $\nu_N$ on the unit disk $\DD$
provided that $2\Re z+2\Re w+N>0$. 

Assume that this condition is satisfied and let $v_{z,w\mid N}$
denote the above function viewed as a vector of the Hilbert space 
$V_N=L^2(\DD,\nu_N)$. Then we have 
$$
\Vert v_{z,w\mid N} \Vert^2=
\int_\DD |(1+D)^z(1+\overline D)^w|^2\nu_N(dD)=
\frac{\Ga(N)\Ga(N+z+\bar z+w+\bar w)}
{\Ga(N+z+\bar w)\Ga(N+\bar z+w)}\,. 
$$
\endproclaim

\demo{Proof} The measure $\nu_N$ is given by the formula of Lemma
3.5. Using it we get
$$
\int_\DD |(1+D)^z(1+\overline D)^w|^2\nu_N(dD)
=\const\, \int_\DD (1+D)^{z+\bar w}(1+\overline D)^{\bar z+w}
(1-|D|^2)^{N-2}\,\ell(dD),\
$$
where $\ell$ denotes the Lebesgue measure. A way of calculating the
latter integral is given in \cite{Ner3, \S1.8}. The constant is fixed
by the requirement that the whole expression equals 1 for $z=w=0$. 
\qed
\enddemo

\proclaim{Lemma 3.14} Assume that $\Re z+\Re w>-\frac12$. Then the
function $f_{z,w\mid N}$ belongs to the Hilbert space
$H_N=L^2(U(N),\mu_N)$, and 
$$
\Vert f_{z,w\mid N}\Vert^2=
\prod_{k=1}^N\frac{\Ga(k)\Ga(k+z+\bar z+w+\bar w)}
{\Ga(k+z+\bar w)\Ga(k+\bar z+w)}\,. 
$$
\endproclaim

\demo{Proof} Using the identification of two spaces indicated in
Lemma 3.7  we can rewrite formula \tht{3.1} as
$$
f_{z,w\mid N}=f_{z,w\mid N-1}\otimes v_{z,w\mid N}\,, \qquad N\ge2.
$$
Note that our assumption on the parameters implies that 
$2\Re z+2\Re w+N>0$ for any $N=2,3,\dots$ (even for $N=1$).
Consequently, $v_{z,w\mid N}$ is well defined as a vector of $V_N$.
Applying Lemma 3.13 we get that $f_{z,w\mid N}$ is square integrable
provided that $f_{z,w\mid N-1}$ is square integrable, and we have by
recurrence 
$$
\Vert f_{z,w\mid N}\Vert=\Vert f_{z,w\mid 1}\Vert\cdot
\Vert v_{z,w\mid 2}\Vert\dots\Vert v_{z,w\mid N}\Vert.
$$

It remains to check that the function $f_{z,w\mid 1}$ is square 
integrable on $U(1)$ (the unit circle) provided that $\Re z+\Re
w>-\frac12$, and 
$$
\Vert f_{z,w\mid 1}\Vert^2=
\frac{\Ga(1)\Ga(1+z+\bar z+w+\bar w)}
{\Ga(1+z+\bar w)\Ga(1+\bar z+w)}\,. 
$$
That is, denoting by $|du|$ the Lebesgue measure on the unit circle
$|u|=1$, 
$$
\int_{|u|=1}(1+u)^{z+\bar w}(1+\bar u)^{\bar z+w}\,\frac{|du|}{2\pi}
=\frac{\Ga(1)\Ga(1+z+\bar z+w+\bar w)}
{\Ga(1+z+\bar w)\Ga(1+\bar z+w)}\,. 
$$
This can be proved by the same argument as in \cite{Ner3, \S1.8}.

Note also that the value of the above integral coincides with the value
of the expression of Lemma 3.13 formally specialized at $N=1$.
This coincidence has an explanation. Indeed, consider the expression
for the measure $\nu_N$ on the disk $\DD$ given in Lemma 3.5, and 
assume that the parameter $N$ is a real number. Then, as $N$ tends to
1 from above, the measure degenerates to the uniform measure
concentrated on the unit circle. 
\qed
\enddemo

The above results make it possible to give the following definition. 

\example{Definition 3.15}

$\bullet$ Let us fix arbitrary $z,w\in\C$. If $N$ is large enough 
($N>-2\Re z-2\Re w$) then the vector $v_{z,w\mid N}\in V_N$ is
well defined, and we can normalize it by setting
$$
v'_{z,w\mid N}=\frac{v_{z,w\mid N}}{\Vert v_{z,w\mid N}\Vert}\,,
$$
where the norm is given by Lemma 3.13. 

$\bullet$ Define an isometric embedding $L_{z,w\mid N}:H_{N-1}\to H_N$
using the identification of Lemma 3.7: 
$$
H_{N-1}\ni f\mapsto f\otimes v'_{z,w\mid N}
\in H_{N-1}\otimes V_N\simeq H^0_N\subset H_N\,.  \tag3.2
$$
It commutes with the action of $G(N-1)$. Consequently, we can form
the inductive limit Hilbert space $H=\varinjlim H_N$ and the natural
unitary representation $\varinjlim\Reg_N$ in $H$. We denote it by
$T_{z,w}$. 

$\bullet$ If $\Re z+\Re w>-\frac12$ then the functions $f_{z,w\mid N}$ are
square integrable and we can normalize them:
$$
f'_{z,w\mid N}=\frac{f_{z,w\mid N}}{\Vert f_{z,w\mid N}\Vert}\,.
$$
Then \tht{3.2} implies existence of the inductive limit vector
$$
\xi_{z,w}=\varinjlim f'_{z,w\mid N}\in H. 
$$
It is $K$--invariant, because all the functions $f_{z,w\mid N}$ are
constant on conjugacy classes. 

$\bullet$ It is worth noting that the definition of this distinguished vector
makes sense only when $\Re z+\Re w>-\frac12$. 
\endexample

\example{Remark 3.16} Removing the normalization of the function
$v_{z,w\mid N}$ we get a map $\wt L_{z,w\mid N}$, which differs from
$L_{z,w\mid N}$ by a scalar multiple and has the form
$$
(\wt L_{z,w\mid N}f)(U)=f(p_N(U))\,v_{z,w\mid N}(U_{NN}), 
\qquad U\in U(N). 
$$
It is worth noting that $\wt L_{z,w\mid N}$ makes sense for any $N$,
irrespective of the values of the parameters $z,w$, and can be
applied to any function $f$ on $U(N-1)$. By virtue of \tht{3.1}, we have
$$
\wt L_{z,w\mid N}: f_{z,w\mid N-1}\mapsto f_{z,w\mid N}\,. 
$$
\qed
\endexample

The formula \tht{3.2} describes the embedding 
$L_{z,w\mid N}: H_{N-1}\to H_N$. More generally, for $M<N$ (where $M$
is large enough), we shall now describe the embedding
$$
L_{z,w\mid M,N}: H_M\to H_N, \qquad
L_{z,w\mid M,N}=L_{z,w\mid N}\circ\dots\circ L_{z,w\mid M+1}\,. \tag3.3
$$

Write any $U\in U(N)$ in block form 
$\bmatrix A & B\\ C & D\endbmatrix$ according to the partition
$N=M+(N-M)$. Then $D$ lies in the matrix ball $\DD_{N-M}$, the set of
complex matrices of size $(N-M)\times(N-M)$, with norm $\le1$.
Generalizing the definitions of $p_N$ and $\ep_N$, consider the
maps 
$$
\gather 
\ep_{M,N}:U(N)\to\DD_{N-M}\,, \qquad U\mapsto D, \\
\wt p_{M,N}: U(N)\to U(M)\times\DD_{N-M}\,, \qquad
\wt p_{M,N}=p_{M,N}\times\ep_{M,N}\,,
\endgather
$$
where $p_{M,N}$ was introduced in Remark 3.2 (if $U\in U(N)'$ then
$p_{M,N}(U)=A-B(1+D)^{-1}C$). Let $\nu_{M,N}$ be the 
image of the Haar measure $\mu_N$ under the map $\ep_{M,N}:
U(N)\to\DD_{N-M}$, and let $V_{M,N}$ denote the Hilbert space
$L^2(\DD_{N-M},\nu_{M,N})$. Next, let 
$H^0_{M,N}\subset H_N=L^2(U(N),\mu_N)$ be the 
subspace of functions depending on $\wt p_{M,N}(\,\cdot\,)$ only. 
As in Lemma 3.7, we have a natural isometry between $H^0_{M,N}$ and
$H_M\otimes V_{M,N}$; we use it to identify these two spaces. 

Set $\DD'_{N-M}=\DD_{N-M}\cap \operatorname{Mat}(N-M)_{>-1}$; this
subset of $\DD_{N-M}$ contains the interior of the matrix ball.
Generalizing the definition of the vector $v_{z,w\mid N}$ (see Lemma 
3.13), we introduce a function $v_{z,w\mid M,N}$ on $\DD_{N-M}$ as
follows: 
$$
v_{z,w\mid M,N}(D)=\cases \det((1+D)^z)\det((1+\overline D)^w), 
& D\in \DD'_{N-M}\,;\\ 0, & D\in\DD_{N-M}\setminus\DD'_{N-M}\,. 
\endcases \tag3.4
$$

\proclaim{Proposition 3.17} Let $M$ be so large that the
embeddings $L_{z,w\mid N}$ are well defined for all $N>M$. Fix $N>M$.

{\rm(i)} The function $v_{z,w\mid M,N}$ defined in \tht{3.4} lies in the
Hilbert space $V_{M,N}=L^2(\DD_{M,N},\nu_{M,N})$.

{\rm(ii)} The embedding $L_{z,w\mid M,N}$ defined in \tht{3.3} maps
$H_M$ into the subspace $H^0_{M,N}$ of $H_N$. Under the
identification $H^0_{M,N}=H_M\otimes V_{M,N}$, this embedding has the
form 
$$
L_{z,w\mid N}: f\mapsto f\otimes v'_{z,w\mid M,N}\,, \qquad
v'_{z,w\mid M,N}=\frac{v_{z,w\mid M,N}}
{\Vert v_{z,w\mid M,N}\Vert}.
$$
\endproclaim

\demo{Proof} {}From the definition \tht{3.3} it follows that 
$$
(L_{z,w\mid N}f)(U)=f(p_{M,N}(U))\, g(U), \tag3.5
$$
where $g$ is a certain function not depending on $f$. We shall prove
that $g(U)$ is proportional to $v_{z,w\mid M,N}(\ep_{M,N}(U))$. It
will follow that the image of $L_{z,w\mid N}$ lies in $H^0_{z,w\mid
M,N}$ and then all the claims of the proposition will become clear.

It is convenient to pass to the maps defined in Remark 3.16. Then we
get, similarly to \tht{3.3}, a map $\wt L_{z,w\mid M,N}$, which differs
from $L_{z,w\mid M,N}$ by a scalar multiple. Similarly to \tht{3.5}, we
have 
$$
(\wt L_{z,w\mid M,N}f)(U)=f(p_{M,N}(U))\, \wt g(U), \tag3.6
$$
where $\wt g$ is proportional to $g$. According to Remark 3.16, we
may apply formula \tht{3.6} to any function $f$, not necessarily a
square integrable one. So, we may take $f=f_{z,w\mid M}$. By the last
claim of Remark 3.16, $\wt L_{z,w\mid N}$ takes
$f_{z,w\mid M}$ to $f_{z,w\mid N}$, which implies
$$
\wt g(U)=\frac{f_{z,w\mid N}(U)}{f_{z,w\mid M}(p_{M,N}(U))}\,.\tag3.6
$$
On the other hand, the same argument as in the proof of Lemma 3.12
shows that the right--hand side is equal to $v_{z,w\mid M,N}(D)$,
where $D=\ep_{M,N}(U)\in\DD_{N-M}$. 

Thus, $\wt g(U)=v_{z,w\mid M,N}(D)$. Since $\wt g$ is proportional to
$g$, this concludes the proof. \qed
\enddemo

\proclaim{Theorem 3.18} All representations constructed in
Definition 3.15 are admissible in the sense of Definition 2.2.
\endproclaim

\demo{Proof} Assume as above that $M$ is a sufficiently large natural
number, so that the embedding $H_M\subset H$ is well defined. By the
very definition, the subspace $\cup_M H_M$ is dense in $H$.

On the other hand, denote by $H'_M$ the subspace in $H$ formed by all
$K_M$-invariant vectors (recall that $K_M\subset K\simeq U(\infty)$
is the subgroup of matrices of the form $\bmatrix 1_M & 0\\ 0 &
*\endbmatrix$, see Definition 2.2). The subspaces $H'_M$ form an
ascending chain (because the subgroups $K_M$ form a descending chain).
We must prove that $\cup_M H'_M$ is dense in $H$. To do this we shall
prove that $H'_M\supset H_M$.

For any $N>M$, let $K_M(N)$ denote the subgroup in 
$K\simeq U(\infty)$ that is the intersection of $K_M$ with $U(N)$.
Clearly, $K_M(N)$ is isomorphic to $U(N-M)$ and we have
$K_M=\cup_{N>M} K_M(N)$. Note that the function $v_{z,w\mid M,N}$ on
the matrix ball $\DD_{N-M}$ is invariant with respect to conjugation
by unitary matrices from $U(N-M)$. Note also that, for any $U\in
U(N)$, the matrix $p_{M,N}(U)$ does not change when we conjugate $U$
by a matrix from $K_M(N)\subset U(N)$. Combining this with the
description of the embedding $H_M\to H_N$ we conclude that the
vectors from $H_M$ are invariant with respect to $K_M(N)$. Since this
holds for any $N$, it follows that all these vectors are
$K_M$-invariant, i.e., $H_M\subset H'_M$.
\qed
\enddemo

\head 4. The space $\U$ of virtual unitary matrices and another
description of the representations $T_{z,w}$ \endhead

The construction of the representations $T_{zw}$ described in \S3 is
quite similar to the second construction of the ``generalized regular
representations'' $T_z$ in \cite{KOV}. In this section we present the
counterpart of the first construction from \cite{KOV}.

\example{Definition 4.1} Let $\U=\varprojlim U(N)$ be the projective
limit of the spaces $U(N)$, as $N\to\infty$, taken with respect to
the projections $p_N: U(N)\to U(N-1)$. By analogy with the space of
virtual permutations from \cite{KOV} we call $\U$ the {\it space of
virtual unitary matrices.\/} 
\endexample

By definition, a point of $\U$ is an arbitrary sequence $x=(x_N)$,
where $x_N\in U(N)$ for any $N=1,2,\dots$, and $p_N(x_N)=x_{N-1}$ for
any $N\ge2$. We denote by $\pi_N$ the natural projection map 
$\U\to U(N)$ sending $x$ to $x_N$. 

There is a natural embedding $U(\infty)\hookrightarrow\U$ assigning
to a matrix $U\in U(\infty)$ a sequence $x=(x_N)$ such that $x_N=U$
provided that $N$ is so large that $U\in U(N)\subset U(\infty)$.
Thus, the image of $U(\infty)$ in $\U$ consists of the stabilizing
sequences $x=(x_N)$. 

However, general elements of the space $\U$ cannot be interpreted as
unitary matrices. 

Recall that the map $p_N$ is
continuous on the open subset $U(N)'\subset U(N)$ but discontinuous
on the whole space $U(N)$ (Remark 3.3).
This is an obstacle to equipping the space $\U$ 
with a natural topology. However, $p_N$ certainly are Borel maps, so
that $\U$ has a natural Borel structure. We keep in mind this
structure while speaking of measures on $\U$.

\example{Definition 4.2} By Lemma 3.11, $p_N$ maps $U(N)'$ onto
$U(N-1)'$, which makes it possible to define the following subset in
$\U$: 
$$
\U'=\varprojlim U(N)'\subset \U.
$$
\endexample

\example{Definition 4.3} A subset of $\U$ will be called {\it
negligible\/} if for any $N$, its image under $\pi_N:\U\to U(N)$ is
a null set with respect to the Haar measure on $U(N)$. 
\endexample

For instance, $\U\setminus\U'$ is a negligible set. Let us agree 
to identify functions on $\U$ that coincide outside a
negligible set. Let also agree that a function is well defined on
$\U$ if it is defined outside a negligible set.

\example{Definition 4.4} We define a right action $(x,g)\mapsto x.g$
of the group $G=U(\infty)\times U(\infty)$ on the space $\U$ by
$$
\pi_N(x.g)=U_2^{-1}\pi_N(x)U_1, \qquad 
x\in\U, \quad g\in (U_1,U_2)\in G,
$$
where $N$ is so large that $g\in G(N)=U(N)\times U(N)$. The
correctness of this definition follows from the basic equivariance
property of the canonical projection $p_N$, see Lemma 3.1.
\endexample

In other words, this action arises from the right actions of the
groups $G(N)$ on the spaces $U(N)$ defined in the beginning of \S3. On
$U(\infty)\subset\U$, this action coincides with 
$(U,\, (U_1,U_2))\to U_2^{-1}UU_1$. Note
also that the action of the subgroup $K=\diag U(\infty)$ arises from
the actions of the groups $U(N)$ on themselves by conjugations.  

Note that the shift of a negligible set by an element of $G$ is
negligible, too. 

\example{Definition 4.5} A function on $\U$ is called {\it
cylindrical\/} if it has the form $F(x)=F_N(\pi_N(x))$ for a certain $N$
and a certain function $F_N$ on $U(N)$. More generally, we extend (in
an evident way) this definition to functions defined outside a
negligible set. 
\endexample

\proclaim{Lemma 4.6} Let $z,w\in\C$, $g\in G$, and
$x\in\U'\cap(\U'.g^{-1})$. If $N$ is large enough then
the expression 
$$
\frac{f_{z,w\mid N}(\pi_N(x.g))}{f_{z,w\mid N}(\pi_N(x))}
$$
does not depend on $N$.
\endproclaim

\demo{Proof} Indeed, assume $g\in G(N-1)$, and let us prove that the
expression above does not change when $N$ is replaced by $N-1$.
Recall that the function $f_{z,w\mid N}$ is the product of
two determinants (which correspond to the particular cases $w=0$ and
$z=0$, respectively). We shall check the above claim for the first
determinant; for the second determinant it is proved similarly, so
that the desired claim will follow.

Let $X=\pi_N(x)$ and $g=(U_1,U_2)$; then $\pi_N(x.g)=U_2^{-1}XU_1$.
Since $g\in G(N-1)$, the matrices $U_1,U_2$ lie in $U(N-1)$. It
follows that $p_N(U_2^{-1}XU_1)=U_2^{-1}p_N(X)U_1$ (Lemma 3.1) and
that the matrix entry $D=X_{NN}$ does not change when $X$ is 
replaced by $U_2^{-1}XU_1$. 

Then, applying Lemma 3.12, we get that the ratio 
$$
\frac{\det((1+U_2^{-1}XU_1)^z)}{\det((1+X)^z)}
$$
does not change when $X$ is replaced by $p_N(X)$, which concludes the
proof. \qed
\enddemo

As a corollary we get

\proclaim{Proposition 4.7} Fix arbitrary $z,w\in\C$. For $g\in G$ and
$x\in\U$, set 
$$
C_{z,w}(x,g)=\text{\rm the stable value of} \quad 
\frac{f_{z,w\mid N}(\pi_N(x.g))}{f_{z,w\mid N}(\pi_N(x))}
\quad \text{\rm as $N\to\infty$.} 
$$
For any $g$, this is a correctly defined cylindrical function in $x$
in the sense of Definition 4.5.
Furthermore, $C_{z,w}(x,g)$ possesses the multiplier property
$$
C_{z,w}(x,g_1)C_{z,w}(x.g_1, g_2)=C_{z,w}(x,g_1g_2), \qquad x\in\U,
\quad g_1,g_2\in G.
$$
Finally,
$$
C_{z,w}(\,\cdot\,, g)\equiv1 \quad \text{\rm for $g\in K$.}
$$
\endproclaim

\demo{Proof} The first claim follows from Lemma 4.6. The second claim
is evident. The third claim follows from the fact that the function
$f_{z,w\mid N}$ on $U(N)$ is central (constant on conjugacy classes).
\qed 
\enddemo

Fix an arbitrary $s\in\C$ with $\Re s>-\frac12$. For each $N=1,2,\dots$
we consider a measure $\mu^{(s)}_N$ on $U(N)$ such that
$$
\mu^{(s)}_N(dU)=\left(\prod_{k=1}^N\frac{\Ga(k)\Ga(k+s+\bar s)}
{\Ga(k+s)\Ga(k+\bar s)} 
\right)^{-1}
|\det((1+U)^s)|^2 \mu_N(dU), 
$$
where, as above, $\mu_N$ is the normalized Haar measure on $U(N)$.  
When $s=0$, this measure reduces to $\mu_N$. Further, using the
formula of Corollary 3.9 and Lemma 3.14, we get that 
$$
\frac{\mu^{(s)}_N(dU)}{\mu_N(dU)}=\frac{|f_{z,w\mid N}(U)|^2}
{\Vert f_{z,w\mid N}\Vert^2} \qquad 
\text{for any $z,w$ such that $z+\bar w=s$.}
$$
This implies that $\mu^{(s)}_N$ is a probability measure. Since
$f_{z,w\mid N}$ is central, $\mu^{(s)}_N$ is invariant under the action
of $U(N)$ on itself by conjugations.

\proclaim{Lemma 4.8} The family
$\{\mu^{(s)}_N\}_{N\ge1}$ is consistent with the projections $p_N$, i.e.,
the pushforward of $\mu^{(s)}_N$ under $p_N$ is $\mu^{(s)}_{N-1}$. 
\endproclaim

\demo{Proof} Same reasoning, based on Lemma 3.12, as in the proof of
Lemma 3.13. \qed
\enddemo

This makes it possible to give the following

\example{Definition 4.9} For any $s\in\C$, $\Re s>-1/2$, we denote by
$\mu^{(s)}$ the projective limit of the family $\{\mu^{(s)}_N\}$ as
$N\to\infty$. This is a probability Borel measure on $\U$.
\endexample

Note the following properties of the measures $\mu^{(s)}$:

$\bullet$ $\mu^{(s)}$ is invariant with respect to the action of $K$ on
$\U$; 

$\bullet$ the measure $\mu^0$ is the projective limit of the Haar
measures $\mu_N$; consequently, it is $G$--invariant;

$\bullet$ negligible sets are null sets with respect to any measure
$\mu^{(s)}$. 

\proclaim{Proposition 4.10} For any $s\in\C$ with $\Re s>-\frac12$, 
the measure $\mu^{(s)}$ is $G$-quasiinvariant.
Its Radon--Nikodym derivative 
$$
\tau_s(x,g)=\frac{\mu^{(s)}(d(x.g))}{\mu^{(s)}(dx)}\,, 
\qquad x\in\U, \quad g\in G,
$$
is a cylindrical function in $x$. Moreover, for any complex $z,w$
such that $z+\bar w=s$ we have 
$$
\tau_s(x,g)=|C_{z,w}(x,g)|^2, \qquad x\in\U, \quad g\in G,
$$
where $C_{z,w}(x,g)$ is the multiplier introduced in Proposition 4.7.
\endproclaim

\demo{Proof} Indeed, this follows from the definition of the measures
$\mu^{(s)}_N$ and Proposition 4.7. \qed
\enddemo

\proclaim{Theorem 4.11} Let $s, z, w\in\C$ be such that $\Re
s>-\frac12$, $z+\bar w=s$, and let $C_{z,w}(x,g)$ be the multiplier
introduced in Proposition 4.7. The following formula defines a unitary
representation $\un T_{z,w}$ of the group $G=U(\infty)\times U(\infty)$ in
the Hilbert space $\un H=L^2(\U,\mu^{(s)})$:
$$
(\un T_{z,w}(g)f)(x)=f(x.g)C_{z,w}(x,g), 
\qquad g\in G, \quad f\in \un H, \quad x\in\U.
$$

Further, the constant function ${\bold 1}\in L^2(\U,\mu^{(s)})$ is a
unit $K$--invariant vector in $\un H$.
\endproclaim

We call ${\bold 1}$ the {\it distinguished vector\/} of
the representation $\un T_{z,w}$.  

\demo{Proof} The correctness of the definition follows from
Proposition 4.10, and the second claim follows from the
$K$--invariance property of the multiplier, see the last claim of
Proposition 4.7. \qed
\enddemo

\proclaim{Proposition 4.12} Let $z,w\in\C$ and $s=z+\bar w$. Assume that
$\Re s>-\frac12$. There exists a Hilbert space isometry between the 
space $\un H$ of the representation $\un T_{z,w}$ and the space $H$ of the
representation $T_{z,w}$ constructed in section 3; this isometry
intertwines the representations $\un T_{z,w}$ 
and $T_{z,w}$ and takes the distinguished vector ${\bold 1}\in \un H$ to the
distinguished vector $\xi_{z,w}\in H$ introduced in Definition 3.15. 
\endproclaim

\demo{Proof} Recall that $H=\varinjlim H_N$, where
$H_N=L^2(U(N),\mu_N)$. On the other hand, since
$\mu^{(s)}=\varprojlim\mu^{(s)}_N$, we have $\un H=\varinjlim \un H_N$, where
$\un H_N=L^2(U(N),\mu^{(s)}_N)$. For each $N=1,2,\dots$, the operator of
multiplication by the function 
$$
f'_{z,w\mid N}(\,\cdot\,)=\frac{f_{z,w\mid N}(\,\cdot\,)}
{\Vert f_{z,w\mid N}\Vert}
$$ 
defines an isometry $\un H_N\to H_N$ taking the constant function 1
to the vector $f'_{z,w\mid N}$. By the very construction, these
isometries are consistent with the embeddings 
$\un H_N\to \un H_{N+1}$ and $H_N\to H_{N+1}$ and, consequently,
define an isometry $\un H\to H$. Clearly, 
the isometry $\un H\to H$ takes ${\bold 1}$ to 
$\xi_{z,w}$. Moreover, again by the construction of Theorem 4.11,
$\un H\to H$ is an intertwining operator between the representations
$\un T_{z,w}$ and $T_{z,w}$. \qed 
\enddemo

\example{Remark 4.13} The above results can be extended, with
appropriate modifications, to the case
when $z,w$ are arbitrary complex numbers (and $s=z+\bar w$, as
usual). When $\Re s\le-\frac12$, the 
measure $\mu^{(s)}$ is infinite and the constant function $\bold 1$
is not a vector of $\un H$. Recall that the
definition of the distinguished vector $\xi_{z,w}$ also fails in this
case. The equivalence of $T_{z,w}$ and $\un T_{z,w}$ remains valid
for all values of the parameters. 
\endexample

\head 5. Measures on the space of Hermitian matrices \endhead

Here we give another description of the measures $\mu^{(s)}$
which were introduced in \S4.  

Let $\Herm(N)$ denote the space of $N\times N$ complex Hermitian matrices.

\example{Definition 5.1} We introduce a bijection $U\leftrightarrow X$
between $U(N)'\subset U(N)$ and $\Herm(N)$ as follows:
$$
U=\frac{i-X}{i+X}\,, \qquad X=i\,\frac{1-U}{1+U}\,,
$$ 
where $i=\sqrt{-1}$ is identified with the scalar matrix 
$i\cdot1_N$. We call $U(N)'\to \Herm(N)$ the {\it Cayley transform\/} and 
$\Herm(N)\to U(N)'$ the {\it inverse Cayley transform.}
\endexample

We define the projection $p'_N: \Herm(N)\to \Herm(N-1)$ as the
operation of removing the $N$th row and the $N$th column from a
$N\times N$ matrix.  

\proclaim{Lemma 5.2} The following diagram is commutative
$$
\CD
U(N)' @>\text{Cayley}>> \Herm(N)\\
@V{p_N}VV @VV{p'_N}V\\
U(N-1)' @>\text{Cayley}>> \Herm(N-1)
\endCD
$$
\endproclaim

\demo{Proof} Indeed, this is verified by a direct calculation. The
easiest way is to use the interpretation of $p_N$ in terms of graphs
of operators, see the proof of Lemma 3.2. \qed
\enddemo

\proclaim{Proposition 5.3} The Cayley transform takes the Haar measure
$\mu_N$ on $U(N)'$ to the measure
$$
\m_N(dX)=\const \det(1+XX^*)^{-N}\Leb(dX)
$$
on $\Herm(N)$. Here $\Leb$ denotes the Lebesgue measure.

More generally, for any $s\in\C$ with $\Re s>-\frac12$, the measure
$\mu^{(s)}_N$ is transformed to the measure 
$$
\m^{(s)}_N(dX)=\const\,\det((1-iX)^{-s})\det((1+iX)^{-\bar s}) 
\det(1+XX^*)^{-N}\Leb(dX).
$$
\endproclaim 

\demo{Proof} Direct computation. \qed
\enddemo

Note that for real values of $s$ the latter expression can be
simplified: 
$$
\m^{(s)}_N(dX)=\const\, \det(1+X^2)^{-s-N}\Leb(dX), \qquad s\in{\Bbb R}.
$$

\proclaim{Corollary 5.4} Fix $s\in\C$, $\Re s>-\frac12$. The measures
$\m^{(s)}_N$ with varying $N$ are consistent with the projections $p'_N$.
\endproclaim

\demo{Proof} Indeed, this follows from Lemma 5.2, Proposition 5.3,
and the similar claim for the measures $\mu^{(s)}_N$. \qed
\enddemo

Of course, a direct verification of Corollary 5.4 is also possible.
For real values of $s$ it was first carried out by Hua, see
the proof of Theorem 2.1.5 in his remarkable book \cite{Hua}. 

The considerations above lead to the following

\example{Definition 5.5} Let $\HH$ denote the space of all infinite
Hermitian matrices. We may view it as a projective limit space:
$\HH=\varprojlim \Herm(N)$. For any $s\in\C$ with $\Re s>-\frac12$, the
family $\{\m^{(s)}_N\}_{N=1,2,\dots}$ defines a probability measure
on the space $\HH$, which will be denoted by $m^{(s)}$ and called the
{\it Hua--Pickrell measure\/} {\rm(}with parameter $s${\rm)}. 
\endexample

The Hua--Pickrell measures are studied in detail in \cite{BO3}. The
measures with real parameter $s$ are exact counterparts of
measures introduced by Pickrell \cite{Pi} (instead of $\HH$, he
considered the space of {\it all\/} infinite complex matrices). 

The group $U(\infty)$ operates on the space $\HH$ by conjugations, and
all the measures $m^{(s)}$ are clearly invariant with respect to this
action. 

The Cayley transforms $U(N)'\to \Herm(N)$ (Definition 5.1) define a
bijection $\U'\to\HH$, which is an isomorphism of measures spaces
$(\U',\mu^{(s)})\to(\HH,m^{(s)})$ for any $s$. Since $\U\setminus\U'$ is a
negligible set, this can also be viewed as an isomorphism 
$(\U,\mu^{(s)})\to(\HH,m^{(s)})$. 

\example{Remark 5.6} The bijection $\U'\to\HH$ takes the action of
$G$ on $\U$ to an action by fractional--linear transformations on
$\HH$. The transformation of a matrix $X\in\HH$ by an element
$g=(U_1,U_2)\in G$ has the form
$$
X\mapsto (Xb+d)^{-1}(Xa+c), \qquad
\bmatrix a & b\\ c & d \endbmatrix=
\bmatrix \frac{U_1+U_2}2 & -i\, \frac{U_1-U_2}2 \\
i\,\frac{U_1-U_2}2 & \frac{U_1+U_2}2 \endbmatrix\,.
$$
Given $g\in G$, this transformation is defined almost everywhere.
\endexample

\example{Remark 5.7} The above results can be extended to the case
of an arbitrary complex parameter $s$, cf. Remark 4.13. However, when
$\Re s\le-1/2$, the measures become infinite. 
\endexample

\head 6. The characters $\chi_{z,w}$ \endhead

Let $\chi$ be a character of $U(\infty)$. Then, for each
$N=1,2,\dots$, the restriction of $\chi$ to $U(N)\subset U(\infty)$
is a character of $U(N)$ (in the sense of Definition 1.1). According
to Example 1.2, we have
$$
\chi\mid_{U(N)}=\sum_{\la\in U(N)\dual} P_N(\la)\wt\chi^\la, \tag6.1
$$
where $P_N(\la)$ are nonnegative coefficients whose sum is equal to 1, and
$\wt\chi^\la$ are the normalized irreducible characters of $U(N)$, see
\tht{1.1}. Note that the series converges uniformly on $U(N)$, because
$|\wt\chi^\la(\,\cdot\,)|\le 1$. 

Any character is uniquely determined by the collection of 
its coefficients $P_N(\la)$. Indeed, the coefficients with fixed $N$ 
determine the restriction of the character to $U(N)$, and if we know
these restrictions for all $N$ then we know the character itself.  

\example{Definition 6.1} Following the general scheme described in
\S2, we introduce a family of characters attached
to the representations $T_{z,w}\,$, as follows. Let $z,w\in\C$ satisfy the
condition $\Re z+\Re w>-\frac12$ ensuring the existence of the distinguished
$K$-invariant vector $\xi_{z,w}$ in the representation $T_{z,w}$, see
Definition 3.15. We consider the matrix coefficient determined by this vector
and then pass to the corresponding character, which we denote by
$\chi_{z,w}$: 
$$
\chi_{z,w}(U)=(T_{z,w}(U,1)\xi_{z,w},\xi_{z,w}), 
\qquad U\in U(\infty).
$$
\endexample

\proclaim{Lemma 6.2} Assume $U\in U(N)\subset U(\infty)$. Then 
$$
\chi_{z,w}(U)=\frac1{\Vert f_{z,w\mid N}\Vert^2}\,
\int_{U(N)}f_{z,w\mid N}(VU)\overline{f_{z,w\mid N}(V)}\mu_N(dV).
$$
\endproclaim

\demo{Proof} By the very construction of the representations
$T_{z,w}$, if $U\in U(N)$ then
$$
\chi_{z,w}(U)=(\Reg_N(U,1)f'_{z,w\mid N}, f'_{z,w\mid N})
=\int_{U(N)}f'_{z,w\mid N}(VU)\overline{f'_{z,w\mid N}(V)}\mu_N(dV),
$$
which is equivalent to the desired formula. \qed
\enddemo

Our aim is to describe explicitly the expansion \tht{6.1} of the characters
$\chi_{z,w}$. 

We shall interpret the labels $\la\in U(N)\dual$ of irreducible
characters of $U(N)$ as {\it signatures of length\/} $N$, i.e., as
ordered $N$-tuples of integers:
$$
\la=(\la_1\ge\dots\ge\la_N), \qquad \la_i\in\Z.
$$

\proclaim{Theorem 6.3} Assume $\Re z+\Re w>-\frac12$. Denote by
$P_{z,w\mid N}(\la)$ the coefficients 
in the expansion \tht{6.1} of the characters $\chi_{z,w}$. We have 
$$
P_N(\la\mid z,w)=
(S_N(z,w))^{-1}\,\cdot\,\prod_{i=1}^N 
\left|\frac1{\Gamma(z-\la_i+i)\Gamma(w+N+1+\la_i-i)}\right|^2\,
\cdot\,\Dim^2_N(\la), \tag6.2
$$
where
$$
S_N(z,w)=\prod_{i=1}^N\frac{\Gamma(z+\bar z+w+\bar w+i)}
{\Gamma(z+w+i)\Gamma(z+\bar w+i)
\Gamma(\bar z+w+i)\Gamma(\bar z+\bar w+i)\Gamma(i)}\,, \tag6.3
$$
and $\Dim_N(\la)$ is the dimension of the irreducible character
$\chi^\la$, given by Weyl's formula:
$$
\Dim_N(\la)=\prod_{1\le i<j\le N}\frac{\la_i-\la_j+j-i}{j-i}\,.
$$
\endproclaim

The proof of Theorem 6.3 is partitioned into a few lemmas. 

Note that the irreducible characters $\chi^\la$, where $\la$ ranges
over $U(N)\dual$, form an orthonormal basis in the subspace of
$H_N=L^2(U(N),\mu_N)$ constituted by central functions. Therefore,
there exists an expansion
$$
f_{z,w\mid N}=\sum_{\la\in U(N)\dual}c_{z,w\mid N}\,\chi^\la
$$
with certain coefficients $c_{z,w\mid N}$.

\proclaim{Lemma 6.4}  We have
$$
P_N(\la\mid z,w)=\frac{|c_{z,w\mid N}|^2}
{\Vert f_{z,w\mid N}\Vert^2}\,.
$$
\endproclaim

\demo{Proof} This follows from Lemma 6.2 and the formula
$$
\int_{U(N)}\chi^\la(VU)\overline{\chi^\mu(V)}\mu_N(dV)
=\delta_{\mu\nu}\,\wt\chi^\la(U), \qquad \la,\mu\in\GT_N\,,
$$
which in turn follows from Schur's orthogonality relations. \qed
\enddemo

Since the norm $\Vert f_{z,w\mid N}\Vert$ is known (Lemma 3.14), our
problem is entirely reduced to evaluating the coefficients
$c_{z,w\mid N}$. 

\proclaim{Lemma 6.5} Assume that
$$
f_1(u)=\sum_{l=-\infty}^\infty c(l)u^l, \qquad c(l)\in\C,
$$
is a function on the unit circle $|u|=1$. For any $N=1,2,\dots$
consider the following function on $U(N)$, which is constant on
conjugacy classes:
$$
f_N(U)=f_1(u_1)\dots f_1(u_N),
$$
where $U$ ranges over $U(N)$ and $(u_1,\dots,u_N)$ denotes the 
spectrum of $U$. Expand $f_N$ in the ``Fourier series'' on the
irreducible characters, 
$$
f_N=\sum_{\la\in U(N)\dual}c(\la)\chi^\la\,.
$$
Then the ``Fourier coefficients'' $c(\la)$ are given by determinants
of order $N$, 
$$
c(\la)=\det[c(\la_i-i+j)]_{1\le i,j\le N}.
$$
\endproclaim

\demo{Proof} This is a well--known combinatorial fact, based on
Weyl's character formula. See, e.g., \cite{Vo, Lemme 2} or \cite{Hua,
Theorem 1.2.1}. \qed
\enddemo

Observe that, by Definition 3.10, 
$$
f_{z,w\mid N}(\diag(u_1,\dots,u_N))
=f_{z,w\mid 1}(u_1)\dots f_{z,w\mid 1}(u_N),
$$
as in Lemma 6.5. By virtue of this lemma, the 
problem of evaluating the ``Fourier coefficients'' of $f_{z,w\mid N}$
is split into two parts: first, calculate the Fourier coefficients of
$f_{z,w\mid 1}$ and, second, calculate the above determinants.  

\proclaim{Lemma 6.6} The coefficients of the Fourier expansion
$$
f_{z,w\mid 1}(u)=\sum_{l=-\infty}^\infty c_{z,w\mid 1}(l)u^l
$$
have the form
$$
c_{z,w\mid 1}(l)=\frac{\Ga(1+z+w)}{\Ga(1+z-l)\Ga(1+w+l)}\,.
$$
\endproclaim

\demo{Proof} Setting $u=e^{i\theta}$ we have
$$
c_{z,w\mid 1}(l)=\frac1{2\pi}\,\int_{-\pi}^\pi
(1+u)^z(1+\bar u)^w u^{-l}d\theta.
$$
To evaluate this integral one can, e.g., reduce it to \cite{Er,
1.5 (30)}. See also \cite{Ner3}. \qed
\enddemo 

\proclaim{Lemma 6.7} The determinants of Lemma 6.5 corresponding to
the Fourier coefficients of Lemma 6.6 can be
explicitly calculated:
$$
c_{z,w\mid N}(\la)=\prod_{i=1}^N
\frac{\Ga(z+w+i)\Ga(i)}{\Ga(z-\la_i+i)\Ga(w+N+1+\la_i-i)}\cdot\Dim_N(\la).
\tag6.4
$$
\endproclaim

\demo{Proof} We will reduce the determinant in question to a known
one, given in \cite{Kra, Lemma 3}: 
$$
\det\left[\prod_{k=1}^{j-1}(x_i+a_k)\cdot
\prod_{k=j}^{N-1}(x_i+b_k)\right]_{i,j=1}^N
=\prod_{1\le i<j\le N}(x_i-x_j)\cdot
\prod_{1\le i\le j\le N-1}(a_i-b_j). \tag6.5
$$

To do this we abbreviate $l_i=\la_i-i$ and we transform our
determinant as follows:
$$
c_{z,w\mid N}(\la)=\frac{(\Ga(1+z+w))^N}
{\prod_{i=1}^N\Ga(z-\la_i+i)\Ga(w+N+1+\la_i-i)}
\cdot\det[A(i,j)],
$$
where
$$
\gather
A(i,j)=\frac{\Ga(z-l_i)\Ga(w+N+1+l_i)}
{\Ga(1+z-l_i-j)\Ga(w+1+l_i+j)}\\
=\underbrace{(z-l_i-1)\dots(z-l_i-j+1)}_{j-1}\cdot
\underbrace{(w+N+l_i)\dots(w+l_i+j+1)}_{N-j}.
\endgather
$$
Then the reduction to \tht{6.5} is carried out by setting 
$$
x_i=l_i, \quad a_j=-z+j, \quad b_j=w+1+j.
$$
Then, after simple transformations, we get \tht{6.4}. \qed

\enddemo

\demo{Proof of Theorem 6.3} Follows from the lemmas above and the formula
for the norm in Lemma 3.14. \qed
\enddemo

\example{Remark 6.8} Note some symmetry properties of the characters
$\chi_{z,w}$.

$\bullet$ Conjugation:
$$
\overline{\chi_{z,w}(U)}=\chi_{\bar w, \bar z}(U), 
\qquad U\in U(\infty). 
$$
Indeed, this follows from Lemma 6.2 and the fact that
$$
\overline{f_{z,w\mid N}(U)}=f_{\bar w, \bar z\mid N}(U), 
\qquad U\in U(N). 
$$

$\bullet$ The character $\chi_{z,w}$ is invariant
under the symmetries of the parameters generated by $z\to\bar z$
and $w\to\bar w$.
Indeed, this follows at once from the corresponding
symmetry of the coefficients $P_N(\la\mid z,w)$. 
Note that this fact is not obvious from the construction of the
characters $\chi_{z,w}$. 

$\bullet$ Let $\la^*=(-\la_N,\dots,-\la_1)$ denote the
{\it dual signature\/} to $\la$: this is the label of the conjugate
irreducible character $\overline{\chi^\la}$. Another symmetry property of the
coefficients $P_N(\la\mid z,w)$ is as follows:
$$
P_N(\la\mid z,w)=P_N(\la^*\mid w,z).
$$
\endexample

Theorem 6.3 is an analog of Proposition 4.8 in \cite{Pi}. This
theorem leads to an important consequence for the representations
$T_{z,w}$, cf. Lemma 4.5 in \cite{Pi}. 

\proclaim{Corollary 6.9} Let $z,w\in\C$ be such that 
$\Re z+ \Re w>-1/2$ and $z,w\notin\Z$. Then the distinguished vector
$\xi_{z,w}$ {\rm(}see Definition {\rm3.15)} is cyclic.
\endproclaim

\demo{Proof} Indeed, for nonintegral values of $z,w$, the coefficients
$P_N(\la\mid z,w)$ are nonvanishing for all $N$ and all
$\la\in U(N)\dual$. This implies that the vector $f_{z,w\mid N}\in H_N$ is
a cyclic vector of the representation $\Reg_N$ for any $N$. This
concludes the proof. \qed
\enddemo

\example{Remark 6.10} Another possible way of evaluating the
coefficients $c_{z,w\mid N}$ is to use Neretin's results, see
\cite{Ner4}. 
\endexample

\head 7. Other series of characters
\endhead 

For two signatures $\nu$ and $\la$, of length $N-1$ and $N$,
respectively, write $\nu\prec\la$ if
$$
\la_1\ge\nu_1\ge\la_2\ge\nu_2\ge\dots\ge\nu_{N-1}\ge\la_N\,.
$$  
The relation $\nu\prec\la$ appears in the {\it Gelfand--Tsetlin
branching rule\/} for the irreducible characters of the unitary
groups, see, e.g., \cite{Zh}: 
$$
\chi^\la\mid_{U(N-1)}=\sum_{\nu:\,\nu\prec\la}\chi^\nu\,.
$$

\example{Definition 7.1}
The {\it Gelfand--Tsetlin graph\/} $\GT$ is a $\Z_+$--graded graph
whose $N$th level $\GT_N$ consists of signatures of length $N$.
Two vertices $\nu\in\GT_{N-1}$ and $\la\in\GT_N$ are connected by an edge if 
$\nu\prec\la$. We agree that $\GT_0$ consists of a single element
denoted as $\varnothing$; it is connected to all $\la\in\GT_1$. \qed
\endexample

\example{Definition 7.2} For $\nu\in\GT_{N-1}$ and $\la\in\GT_N$,
where $N=1,2,\dots$, set
$$
q(\nu,\la)=\cases 
\dfrac{\Dim_{N-1}\nu}{\Dim_N\la}\,, & \nu\prec\la,\\
0, & \nu\nprec\la. \endcases
$$
This is the {\it cotransition probability function\/} of the
Gelfand--Tsetlin graph; it satisfies the relation 
$$
\sum_{\nu\in\GT_{N-1}}q(\nu,\la)=1, \qquad \forall \la\in\GT_N\,.
$$
We agree that $q(\varnothing, \la)=1$ for all $\la\in\GT_1$. \qed
\endexample

One can imagine that each edge $(\nu,\mu)$ of the graph $\GT$ is
equipped with a label, which is the positive number $q(\nu,\la)$. Note
that the Gelfand--Tsetlin branching rule is equivalent to the
relation 
$$
\wt\chi^\la\mid_{U(N-1)}=\sum_{\nu\in\GT_{N-1}}q(\nu,\la)\wt\chi^\nu\,.
\tag7.1
$$

\example{Definition 7.3} Assume that for each $N=0,1,\dots$ we are
given a probability measure 
$P_N$ on the discrete set $\GT_N$. The family $\{P_N\}$ is
called a {\it coherent system\/} if 
$$
P_{N-1}(\nu)=\sum_{\la\in\GT_N}q(\nu,\la)P_N(\la),
\qquad N=1,2,\dots, \quad \nu\in\GT_{N-1} \,.
\tag7.2
$$
\qed
\endexample

Note that if $P_N$ is an arbitrary probability measure on $\GT_N$
then \tht{7.2} defines a probability measure on $\GT_{N-1}$
(indeed, this follows at once from the above relation for
$q(\nu,\la)$). Thus, in a coherent system $\{P_N\}$, the
$N$th term is a refinement of the $(N-1)$th one. 

\proclaim{Proposition 7.4} There exists a natural bijective correspondence
$\chi\leftrightarrow \{P_N\}$ between characters of the group
$U(\infty)$ and coherent systems on the Gelfand--Tsetlin graph. The
correspondence is defined by the relations
$$
\chi_N=\sum_{\la\in\GT_N}P_N(\la)\,\wt\chi^\la, \qquad 
\chi_N:=\chi\mid_{U(N)}, \quad N=1,2,\dots\,. \tag7.3
$$
\endproclaim

\demo{Proof} Let $\chi$ be a character of $U(\infty)$. 
We repeat the argument at the beginning of \S6. Since $\chi_N$ is a
character of $U(N)$ in the sense of Definition 1.1 ($N=1,2,\dots$),
we get, according to Example 1.2, the expansion \tht{7.3}, where
$P_N(\,\cdot\,)$ is a probability measure on 
$\GT_N$. Then the relation \tht{7.2} follows from \tht{7.1} and the evident
relation $\chi_N\mid_{U(N-1)}=\chi_{N-1}$.

Conversely, let $\{P_N\}$ be a coherent system. We define by means of
\tht{7.3} a sequence $\{\chi_N\}$ of characters of the groups $U(N)$. The
coherency property \tht{7.2} ensures that $\chi_N\mid_{U(N-1)}=\chi_{N-1}$
for any $N=2,3,\dots$, so that there exists a function $\chi$ on
$U(\infty)$ such that $\chi\mid_{U(N)}=\chi_N$ for all
$N=1,2,\dots$. Obviously, $\chi$ is a character of $U(\infty)$. \qed
\enddemo

Using Proposition 7.4 we will construct new series of characters
by analytic continuation of the formulas \tht{6.2}--\tht{6.3}.

Let $z,z',w,w'$ be complex parameters. 
For any $N=1,2,\dots$ and any $\la\in\GT_N$ set
$$
\gathered
P'_N(\la\mid z,z'w,w')=\Dim_N^2(\la)\\
\times \prod_{i=1}^N 
\frac1{\Gamma(z-\la_i+i)\Gamma(z'-\la_i+i)
\Gamma(w+N+1+\la_i-i)\Gamma(w'+N+1+\la_i-i)}\,,
\endgathered
$$
where $\Dim_N\la$ was defined in \S6. Clearly, for any fixed $N$ and $\la$,
$P'_N(\la\mid z,z',w,w')$ is an entire function on $\C^4$. This
expression was obtained by analytic continuation from the expression
given in \tht{6.2}. 

Set
$$
\Cal D=\{(z,z',w,w')\in\C^4\mid \Re(z+z'+w+w')>-1\}.
$$
This is a domain (a half--space) in $\C^4$.

\proclaim{Proposition 7.5} Fix an arbitrary $N=1,2,\dots$\,. 
The series of entire functions
$$
\sum_{\la\in\GT_N} P'_N(\la\mid z,z',w,w') \tag7.4
$$
converges in the domain $\Cal D$, uniformly on compact sets. Its sum
is equal to
$$
S_N(z,z',w,w')=\prod_{i=1}^N
\frac{\Gamma(z+z'+w+w'+i)}
{\Gamma(z+w+i)\Gamma(z+w'+i)\Gamma(z'+w+i)\Gamma(z'+w'+i)\Gamma(i)}
\tag7.5 
$$
\endproclaim

\demo{Proof}  Let us prove that for any $\la\in\GT_N$
$$
|P'_N(\la\mid\zw)|
\le \const\,\prod_{i=1}^N(1+|\la_i|)^{-\Re(z+z'+w+w'+2)}\,, \tag7.6
$$
uniformly on compact sets in $\Cal D$. 

Indeed, using the formula $(\Ga(z)\Ga(1-z))^{-1}=\sin(\pi z)/\pi$ we get
$$
\gather
\frac1{\Gamma(z-\la_i+i)\Gamma(z'-\la_i+i)
\Gamma(w+N+1+\la_i-i)\Gamma(w'+N+1+\la_i-i)}\\
=\frac{\sin(\pi z)\sin(\pi z')}{\pi^2}\,
\frac{\Ga(-z+1+\la_i-i)\Ga(-z'+1+\la_i-i)}
{\Gamma(w+N+1+\la_i-i)\Gamma(w'+N+1+\la_i-i)}\\
=\frac{\sin(\pi w)\sin(\pi w')}{\pi^2}\,
\frac{\Ga(-w-N-\la_i+i)\Ga(-w'-N-\la_i+i)}
{\Gamma(z-\la_i+i)\Gamma(z'-\la_i+i)}\,.
\endgather
$$
Suppose $|\la_i|\to\infty$. Then, applying the first or the second
equality above (depending on whether $\la_i\gg0$ or $\la_i\ll0$) and
the asymptotic formula
$$
\lim_{x\to+\infty}\frac{\Ga(a+x)}{\Ga(b+x)}=x^{a-b}\,,
\qquad a,b\in\C,
$$ 
we get  
$$
\gather
\left|\frac1{\Gamma(z-\la_i+i)\Gamma(z'-\la_i+i)
\Gamma(w+N+1+\la_i-i)\Gamma(w'+N+1+\la_i-i)}\right|\\
\le\const\,(1+|\la_i|)^{-\Re(z+z'+w+w'+2N)}\,, \tag7.7
\endgather
$$
where the estimate is uniform on compact sets in $\Cal D$.

Next, Weyl's formula implies that $\Dim^2_N\la$ is a polynomial in
$\la_1,\dots,\la_N$ of degree $2N-2$ with respect to each variable.
Combining this fact with \tht{7.7} we get \tht{7.6}.

The bound \tht{7.6} ensures the convergence of the series \tht{7.4},
uniformly on compact sets in $\Cal D$, to a holomorphic function in
$\Cal D$. It remains to show that this function coincides with
\tht{7.5}. To do this we remark that $\tht{7.5}$ is a holomorphic
function coinciding with the function $S_N(z,w)$ (see \tht{6.3}) on
the subset $z'=\bar z$, $w'=\bar w$ in $\Cal D$. For this subset, the
claim holds by virtue of Theorem 6.3. Clearly, this subset is a set of
uniqueness for holomorphic functions in $\Cal D$. This implies that
the claim holds on the whole domain $\Cal D$.
\qed
\enddemo

Note that in the special case $N=1$, the set $\GT_1$ is simply $\Z$
and the identity
$$
\sum_{\la\in\GT_1}P'_1(\la\mid\zw)=S_1(\zw)
$$
is equivalent to the well--known Dougall's formula, see \cite{AAR,
Chapter 2, Theorem 2.8.2 and Exercise 42(b)}, \cite{Er, \S1.4}). 

\example{Definition 7.6} The set of {\it admissible values\/} of the
parameters $\zw$ is the subset $\Dadm$ of the quadruples
$(\zw)\in\Cal D$ such that:

First, $P'_N(\la\mid\zw)\ge0$ for any $N$ and any $\la\in\GT_N$.

Second, for any $N$, the above inequality is strict at least for two
different $\la$'s. 

A quadruple $(\zw)$ will be called {\it admissible\/} if it belongs
to $\Dadm$.
\endexample

Consider the subdomain
$$
\gather
\Cal D_0=\{(z,z',w,w')\in\Cal D\mid
z+w,\,z+w',\,z'+w,\,z'+w'\ne-1,-2,\dots\}\\
=\{(\zw)\in\Cal D\mid S_N(z,z',w,w')\ne0\}. 
\endgather
$$
For any $(\zw)\in\Cal D_0$ we set
$$
P_N(\la\mid\zw)=\frac{P'_N(\la\mid\zw)}{S_N(\zw)}\,, \qquad 
N=1,2,\dots,\quad \la\in\GT_N. \tag7.8
$$

Note that $\Dadm\subset\Cal D_0$ (indeed, if $(\zw)\in\Dadm$ then
$S_N(\zw)$ is strictly positive), so that formula \tht{7.8} makes
sense for any admissible $(\zw)$. 

\proclaim{Proposition 7.7} For any $(\zw)\in\Cal D_0$, the
expressions \tht{7.8} satisfy the coherency relation \tht{7.2}.
\endproclaim

\demo{Proof} Plug in $P_{N-1}(\nu)=P_{N-1}(\nu\mid\zw)$ and
$P_N(\la)=P_N(\la\mid\zw)$ to \tht{7.2}. First of all, since $0\le
q(\nu,\la)\le1$, the series in the right--hand side of \tht{7.2}
converges, uniformly on compact sets in $\Cal D_0$, by virtue of
Proposition 7.5. Then we apply the same argument as in the proof of
Proposition 7.5. Namely, we remark that the relation holds provided
that $z'=\bar z$, $w'=\bar w$, and then conclude that it must hold
for any $(\zw)\in\Cal D_0$ by analytic continuation. \qed
\enddemo

\proclaim{Corollary 7.8} For any $(\zw)\in\Dadm$, the expressions
\tht{7.8} form a coherent system. Hence, by Proposition 7.4, there
exists a character $\chi_\zw$ corresponding to this coherent system.
\endproclaim

Our aim is to describe the set $\Dadm$ explicitly.

Define the subset $\Cal Z\subset\C^2$ as follows:
$$
\gather
\Cal Z=\Zp\sqcup\Zc\sqcup\Zd, \\
\Zp=\{(z,z')\in\C^2\setminus\R^2\mid z'=\bar z\}, \\
\Zc=\{(z,z')\in\R^2\mid \exists m\in\Z, \, m<z,z<m+1\}, \\
\Zd=\underset{m\in\Z}\to{\sqcup}\Zdm, \\
\Zdm=\{(z,z')\in\R^2\mid z=m, \,z'>m-1, \quad 
\text{or} \quad z'=m,\, z>m-1\}, 
\endgather
$$
where ``princ'', ``compl'', and ``degen'' are abbreviations for
``principal'', ``complementary'', and ``degenerate'', respectively.
The terminology is justified by the following lemma.  

\proclaim{Lemma 7.9} Let $(z,z')\in\C^2$. 

{\rm(i)} The expression
$(\Gamma(z-k+1)\Gamma(z'-k+1))^{-1}$ is nonnegative for all $k\in\Z$ if
and only if $(z,z')\in\Cal Z$. 

{\rm(ii)} If $(z,z')\in\Zp\sqcup\Zc$ then this expression is
strictly positive for all $k\in\Z$.

{\rm(iii)} If $(z,z')\in\Zdm$ then this expression vanishes for
$k=m+1,m+2,\dots$ and is strictly positive for $k=m,m-1,\dots$\,. 
\endproclaim

\demo{Proof} Assume first that both $z$ and $z'$ are nonintegral.
Then the expression $(\Gamma(z-k+1)\Gamma(z'-k+1))^{-1}$ does not vanish
for any $k\i\Z$. Clearly, it is strictly positive for all $k\in\Z$
whenever $(z,z')$ is in $\Zp$ or in $\Zc$. Let us check the inverse
claim. Dividing $\Gamma(z-k+1)\Gamma(z'-k+1)$ by
$\Gamma(z-k)\Gamma(z'-k)$ we see that $(z-k)(z'-k)$ must be strictly
positive for all $k\in\Z$. But this implies that $(z,z')$ belongs
either to $\Zp$ or to $\Zc$.

Thus, we have verified all the claims in the case when both $z$ and
$z'$ are nonintegral. Now we shall do the same when at least one of
them is integral. By virtue of the symmetry $z\leftrightarrow z'$, we
may assume that $z=m\in\Z$ and $z'\ne m-1,m-2,\dots$. Then the
expression $(\Gamma(z-k+1)\Gamma(z'-k+1))^{-1}$ vanishes for 
$k=m+1,m+2,\dots$ and does not vanish for $k=m,m-1,\dots$\,. If $z'$
is real and strictly greater than $m-1$ then the expression is
strictly positive for $k=m,m-1,\dots$, because then both $\Gamma(z-k+1)$ and
$\Gamma(z'-k+1)$ are strictly positive. Conversely, let
$\Gamma(z-k+1)\Gamma(z'-k+1)$ be strictly positive for $k=m,m-1,\dots$\,.
As $\Gamma(z-k+1)=\Gamma(m-k+1)$ is strictly positive,
$\Gamma(z'-k+1)$ must be strictly positive, too ($k=m,m-1,\dots$).
Hence the same holds for the ratio $\Gamma(z'-k+1)/\Gamma(z'-k+2)$.
Therefore, $z'-k+1>0$ for all $k=m,m-1,\dots$, which implies
$z'>m-1$. This concludes the proof. 
\qed
\enddemo

\proclaim{Proposition 7.10} The set $\Dadm$ introduced in Definition
7.6 consists of the quadruples $(\zw)\in\Cal D$ satisfying the
following two conditions:

First, both $(z,z')$ and $(w,w')$ belong to $\Cal Z$.

Second, in the particular case when both $(z,z')$ and $(w,w')$ are in
$\Zd$, an extra condition is added: let $k,l$ be such that
$(z,z')\in\Zdk$ and $(w,w')\in\Zdl$; then we require $k+l\ge1$. 
\endproclaim

\demo{Proof} Let us abbreviate $P'_N(\la)=P'_N(\la\mid\zw)$. 

Assume that the above two conditions on $(\zw)\in\Cal D$
are satisfied and prove that $(\zw)\in \Dadm$. Indeed, examine in
detail the possible cases. If none of the parameters is integral then
claim (ii) of Lemma 7.9 shows that $P'_N(\la)$ is always strictly
positive. If $(z,z')\in\Zdm$ while $w,w'$ are nonintegral then claim
(iii) of Lemma 7.9 (applied to the couple $(z,z')$) together with claim
(ii) (applied to the couple $(w,w')$) show that $P'_N(\la)$ is
strictly positive if $\la_1\le m$ and vanishes otherwise. Likewise,
if $(w,w')\in\Zdm$ (for some $m\in\Z$) and $z,z'$ are nonintegral then
$P'_N(\la)$ is strictly positive when $\la_N\ge-m$ and vanishes
otherwise. Finally, if $(z,z')\in\Zdk$ and $(w,w')\in\Zdl$ (with some
$k,l\in\Z$) then $P'_N(\la)$ is strictly positive if $\la$ satisfies
the inequalities $k\ge\la_1\ge\dots\ge\la_N\ge-l$ and vanishes
otherwise. Since $k>-l$ by the second condition, these inequalities
are satisfied by $\ge2$ different $\la$'s. Hence we conclude that in
all cases $(\zw)\in\Dadm$, as required.

Conversely, assume that $(\zw)\in\Dadm$ and prove that the above
two conditions hold. We will verify that $(z,z')\in\Cal Z$. Then the
similar claim concerning $(w,w')$ will follow by virtue of the
symmetry property
$$
P'_N(\la_1,\dots,\la_N\mid\zw)=
P'_N(-\la_N,\dots,-\la_1\mid w,w',z,z').
$$
As for the second condition,
$k+l\ge1$, it follows from the argument above. 

We need two simple lemmas.

\proclaim{Lemma 7.11} Assume that for a given $N$ and a certain
$\la\in\GT_N$ the following two conditions hold:

$\bullet$ $\la^\downarrow:=(\la_1-1,\la_2,\dots,\la_N)\in\GT_N$, i.e.,
$\la_1>\la_2$ if $N\ge2$, 

$\bullet$ $P'_N(\la)>0$ and $P'(\la^\downarrow)>0$.

Then 
$$
\frac{(z-\la_1+1)(z'-\la_1+1)}{(w+N+\la_1-1)(w'+N+\la_1-1)}>0.
$$
\endproclaim

\demo{Proof} Indeed, the above expression coincides with the ratio
$P'_N(\la)/P'_N(\la^\downarrow)$. \qed
\enddemo

\proclaim{Lemma 7.12} Assume that for a given $N$ and a certain
$\la\in\GT_N$ the following two conditions hold:

$\bullet$ $\la^\uparrow:=(\la_1,\dots,\la_{N-1},\la_N+1)\in\GT_N$, i.e.,
$\la_{N-1}>\la_N$ if $N\ge2$, 

$\bullet$ $P'_N(\la)>0$ and $P'(\la^\uparrow)>0$.

Then 
$$
\frac{(w+\la_N+1)(w'+\la_N+1)}{(z-\la_N+N-1)(z'-\la_N+N-1)}>0.
$$
\endproclaim

\demo{Proof} Indeed, the above expression coincides with the ratio
$P'_N(\la)/P'_N(\la^\uparrow)$. \qed
\enddemo

Now we resume the proof of Proposition 7.10. To prove that
$(z,z')\in\Cal Z$ we will examine in succession four possible cases.

{\it Case 1: all parameters are nonintegral.\/} Then $P'_N(\la)$ is
always nonzero, hence we have $P_N(\la)>0$ for all $N$ and all
$\la\in\GT_N$. Given $N=1,2,\dots$ and $k\in\Z$, choose $\la\in\GT_N$
such that $\la_1=k$ and (if $N\ge2$) $\la_1>\la_2$. Then, applying
Lemma 7.11 we get 
$$
\frac{(z-k+1)(z'-k+1)}{(w+N+k-1)(w'+N+k-1)}>0.
$$
Fix $k$ and let $N\to\infty$. Then 
$$
\frac{(z-k+1)(z'-k+1)}{(w+N+k-1)(w'+N+k-1)}\sim
\frac{(z-k+1)(z'-k+1)}{N^2}>0.
$$
This implies $(z-k+1)(z'-k+1)>0$. Since this inequality holds for any
$k\in\Z$ we conclude that $(z,z')$ belongs either to $\Zp$ or to
$\Zc$. 

{\it Case 2: at least one of the parameters $z,z'$ is integral and at
least one of the parameters $w,w'$ is integral, too.\/} Then, using
the symmetries $z\leftrightarrow z'$ and $w\leftrightarrow w'$, we
may assume, without loss of generality, that
$$
z=k\in\Z, \quad z'\ne k-1,k-2,\dots; \qquad
w=l\in\Z, \quad w'\ne l-1,l-2,\dots\,.
$$
Let $N$ be arbitrary. We have: $P'_N(\la)=0$ whenever $\la_1>k$ or
$\la_N<-l$. Next, if $k\ge\la_1\ge\dots\ge\la_N\ge-l$ then
$P'_N(\la)\ne0$. It follows that $k\ge-l$ and
even more, $k>-l$, 
because two different $\la$'s with $P'_N(\la)>0$ must exist. 

Now we apply Lemma 7.11 to $\la=(k,k-1,\dots,k-1)$. The assumptions of
the lemma are satisfied, and we get the same inequality as in Case 1
above. Moreover, as $z-k+1$ reduces to 1, the same argument as above
shows that $z'-k+1>0$, i.e., $z'>k-1$. Hence, $(z,z')\in\Zdk$.

{\it Case 3: at least one of the parameters $z,z'$ is integral while
both $w$ and $w'$ are nonintegral.\/} We may assume that $z=k\in\Z$
and $z'\ne k-1,k-2,\dots$\,. In this case, $P'_N(\la)$ does not
vanish (and hence is strictly positive) whenever $\la_1\le k$. Then
we apply the argument of Case 2 and get exactly the same conclusion. 

{\it Case 4: both $z,z'$ are nonintegral while at least one of the
parameters $w,w'$ is integral.\/} We may assume that $w=l\in\Z$ and
$w'\ne l,l-1,\dots$\,. Then we have $P'_N(\la)>0$ whenever
$\la_N\ge-l$. We apply first Lemma 7.11, where we take any $\la$ such
that $\la_1=-l+i$ with $i=1,2,\dots$, and $\la_2<\la_1$ (if $N\ge2$).
The same argument as above then gives the inequalities 
$$
(z+l-i+1)(z'+l-i+1)>0, \qquad i=1,2,\dots\,.
$$

Next, we apply Lemma 7.12, where we take any $\la$ such that
$\la_N=-l$ and $\la_{N-1}>\la_N$ (if $N\ge2$). This leads to the
inequality 
$$
\frac{(w-l+1)(w'-l+1)}{(z+l+N-1)(z'+l+N-1)}>0, \qquad N=1,2,\dots\,.
$$
Applying the same trick as above ($N\to\infty$) we see that the
numerator must be strictly positive. But then the denominator must be
strictly positive for any $N$. This results in the inequalities 
$$
(z+l+j)(z'+l+j)>0, \qquad j=N-1=0,1,\dots\,.
$$

Combining these two families of inequalities we get that
$(z+k)(z'+k)>0$ for any $k\in\Z$, which means that $(z,z')$ belongs
either to $\Zp$ or to $\Zc$.

This completes the proof of Proposition 7.10.
\enddemo

\head 8. Topology on the space of extreme characters \endhead 

Recall that there is a 1--1 correspondence between extreme characters
of $U(\infty)$ and points of the set $\Om\subset\R^{4\infty+2}$, see
Theorem 1.3 and Remark 1.6. The aim of this section is to prove that
this correspondence is a homeomorphism with respect to natural
topologies on both spaces, the space of extreme characters and the space
$\Om$. This result is used in the proof of Theorem 9.1 below, it is
also of independent interest.

First, we have to define the topologies in question. Let us start with
the space $\Om$. We equip it with the product topology of the ambient
space $\R^{4\infty+2}$. Since $\R^{4\infty+2}$ is a separable
metrizable space, so is $\Om$. {}From the definition of $\Om$ it
follows that this is a locally compact space. Furthermore, for any
positive constant $c$, the subset of the form
$\{\om\in\Om\mid\de^++\de^-\le c\}$ is compact.

Now let us turn to the space of extreme characters. This is a subset
of the space $\Cal X(U(\infty))$ of all characters. We equip $\Cal
X(U(\infty))$ with the topology of uniform convergence on each of
the compact subgroups $U(N)\subset U(\infty)$, $N=1,2,\dots$\,. Then
we restrict this topology to the subspace of extreme characters. 
It is readily seen that the topological space thus obtained is separable
and metrizable.  

\proclaim{Theorem 8.1} The correspondence $\om\mapsto\chi^{(\om)}$
between the points $\om\in\Om$ and the extreme characters of the group
$U(\infty)$ is a homeomorphism with respect to the topologies defined
above. 
\endproclaim

\demo{Proof} 
Recall that any character $\chi\in\Cal X(U(\infty))$ is
uniquely determined by the coefficients $P_N(\la)$ of the expansions
\tht{7.3}, where $N=1,2,\dots$ and $\la$ ranges over $\GT_N$. It is
readily seen that the topology of $\Cal X(U(\infty))$ just defined
coincides with the topology of simple convergence of these ``Fourier
coefficients''. Indeed, the crucial point here is that the
coefficients are nonnegative and each sum of the form $\sum
P_N(\la)$, where $\la$ ranges over $\GT_N$, equals 1.  

By virtue of the multiplicativity property of formula \tht{1.2}, the
topology on extreme characters is defined by convergence of the
coefficients $P_1(\la)$, where $\la\in\GT_1=\Z$. These are simply the
ordinary Fourier coefficients of the functions on the unit circle in
$\C$, which are given by the expression in curved brackets in
\tht{1.2}, i.e.,
$$
F^{(\om)}(u)=e^{\ga^+(u-1)+\ga^-(u^{-1}-1)}
\prod_{i=1}^\infty\frac{1+\be_i^+(u-1)}{1-\al_i^+(u-1)}
\,\frac{1+\be_i^-(u^{-1}-1)}{1-\al_i^-(u^{-1}-1)}\,, 
\qquad u\in\C, \quad |u|=1. \tag8.1
$$

According to the remarks above, the claim of the theorem is
equivalent to the following one: in the set of functions of the form
\tht{8.1}, the uniform convergence on the unit circle (or, which is
the same, the simple convergence of the Fourier coefficients) is
equivalent to the convergence of the labels $\om$ in the space $\Om$. 
Let us notice once again that the reformulation in terms of the
convergence of the Fourier coefficients is possible because the
functions \tht{8.1} are positive definite and normalized at $u=1$. 

We proceed in a few steps.

\demo{Step 1} Let us prove that the map $\om\mapsto F^{(\om)}$ is
continuous. Rewrite \tht{8.1} in the form 
$$
e^{\de^+(u-1)+\de^-(u^{-1}-1)}
\prod_{i=1}^\infty\frac{(1+\be_i^+(u-1))e^{-\be^+_i(u-1)}}
{(1-\al_i^+(u-1))e^{\al^+_i(u-1)}}
\,\frac{(1+\be_i^-(u^{-1}-1))e^{-\be^-_i(u^{-1}-1)}}
{(1-\al_i^-(u^{-1}-1))e^{\al^-_i(u^{-1}-1)}}\,, \tag8.2
$$
where we used the fact that
$\ga^\pm=\de^\pm-\sum(\al^\pm_i+\be^\pm_i)$. It suffices to show that
the infinite product in \tht{8.2} converges uniformly on $u$ and
$(\al^+,\be^+,\al^-,\be^-)$, where $u$ ranges over the unit circle
and the sum of the parameters $\al^\pm_i$, $\be^\pm_i$ is bounded by
a constant. This is reduced to the following elementary fact: the infinite
product $\prod(1+x_i)\exp(-x_i)$ converges uniformly on any subset of
$\C^\infty$ of the form $\{(x_i)\in\C^\infty\mid \sum|x_i|\le c\}$,
where $c$ is an arbitrary positive constant. 

Note that without exponentials, the convergence is nonuniform. For
this reason, we cannot take $\ga^\pm$ instead of $\de^\pm$. Note also
that $\ga^+$ and $\ga^-$ are not continuous functions of $\om$.
\enddemo

\demo{Step 2} For $\om\in\Om$, set $\Vert\om\Vert=\de^++\de^-$. We
claim that the inverse map $F^{(\om)}\to\om$ is continuous provided
that $\om$ is subject to the restriction $\Vert\om\Vert\le c$, where
$c$ is an arbitrary positive constant. Indeed, this follows from
the result of step 1, because any subset of the form $\Vert\om\Vert\le
c$ is compact. Here we use the fact that a bijective continuous map
of a compact space on a Hausdorff space is a homeomorphism.
\enddemo

\demo{Step 3} For $\om\in\Om$, set
$$
\Vert\om\Vert'=\ga^++\ga^-+\sum\al^+_i(1+\al^+_i)+
\sum\be^+_i(1-\be^+_i)+\sum\al^-_i(1+\al^-_i)+\sum\be^-_i(1-\be^-_i).
$$
Note that all summands are nonnegative and the sums are finite. 
Let $(\om_n)$ be a sequence of points of $\Om$ such
that the corresponding sequence of the functions $F^{(\om_n)}$ is
convergent. We claim that then $\Vert\om_n\Vert'$ remains bounded.

Indeed, recall that any function of the form \tht{8.1} is positive
definite on the unit circle and normalized at 1. Hence, it is the
characteristic function (i.e., Fourier transform) of a probability
measure on $\Z$. Since \tht{8.1} is a real analytic function, the
corresponding measure possesses finite moments of any order. Note
also that the uniform convergence of characteristic functions on the
unit circle is equivalent to the weak convergence of the
corresponding probability measures.

We will need the following

\proclaim{Lemma {\rm \cite{OkOl, Lemma 5.2}}} Let $(M_n)$ be a
sequence of probability measures on 
$\Z$ {\rm(}or, even more generally, on $\R${\rm)} such that each
$M_n$ has finite moment of order 4 and the 
sequence $(M_n)$ weakly converges to a probability measure. Assume
that the second moments of $M_n$'s tend to infinity. Then the fourth
moments grow faster than the squares of the second moments. 
\endproclaim

Actually, we will use a corollary of the lemma.

\proclaim{Corollary} Let $(M_n)$ be a sequence of probability
measures on $\Z$ with finite 4-th moments, weakly convergent to a
probability measure. Assume that for any $n$,
the 4--th moment of $M_n$ is bounded by the square of 
the 2--nd moment times a constant which does not depend on $n$. Then the
2--nd moments are uniformly bounded.
\endproclaim

Consider the functions
$G_n(u)=F^{(\om_n)}(u)\overline{F^{(\om_n)}(u)}$. These are also
positive definite normalized functions on the unit circle, and the
sequence $(G_n)$ is uniformly convergent by the assumption on the
initial functions. Let $M_n$ be the probability measures
corresponding to the functions $G_n$. Then the measures $M_n$ weakly
converge to a probability measure. We will prove that these measures
obey the assumption of the corollary, which will imply the uniform
boundedness of their second moments. This, in turn, will imply that
$\Vert\om_n\Vert'$ remains bounded, as required. \footnote{We have
replaced the initial functions by the squares of their moduli,
because this simplifies the estimation of the moments.}

Let us realize this plan. Set $u=e^{i\theta}$. If $M$ is a
probability measure on $\Z$ and $G(u)=G(e^{i\theta})$ is its
characteristic function, then the 2--nd and 4--th moments of $M$
are equal, within number factors, to the 2--nd and 4--th coefficients
of the Taylor expansion of the function $G(e^{i\theta})$, respectively. 

Given $\om\in\Om$, introduce new variables as follows
$$
\{a_i\}=\{\al^+_j(1+\al^+_j)\}\sqcup\{\al^-_k(1+\al^-_k)\}, \quad
\{b_i\}=\{\be^+_j(1-\be^+_j)\}\sqcup\{\be^-_k(1-\be^-_k)\}, \quad
c=\ga^++\ga^-.
$$
Here the ordering of the variables $a_i$, $b_i$ is unessential. 
In this notation we get, after a simple computation,
$$
\gather
F^{(\om)}(e^{i\theta})\overline{F^{(\om)}(e^{i\theta})}
=e^{2c(\cos\theta-1)}
\prod_{i=1}^\infty\frac{1+2b_i(\cos\theta-1)}{1-2a_i(\cos\theta-1)}\\
=e^{-c\theta^2+\frac1{12}c\theta^4+\dots}\,
\prod_{i=1}^\infty\frac{1-b_i\theta^2+\frac1{12}b_i\theta^4+\dots}
{1+a_i\theta^2-\frac1{12}a_i\theta^4+\dots}\,.
\endgather
$$

In this expression, the coefficient of $\theta^2$ equals
$$
-(c+\sum a_i+\sum b_i)=-\Vert\om\Vert'\tag8.3
$$
and the coefficient of $\theta^4$ equals
$$
\frac1{12}(c+\sum a_i+\sum b_i)+\frac12 c^2+\sum a_i^2+
e_2(c,a_1,a_2,\dots,b_1,b_2,\dots),\tag8.4
$$
where $e_2$ stands for the 2--nd elementary symmetric function. It is
readily seen that the expression \tht{8.4} is bounded from above by
the square of the expression \tht{8.3} times a constant. This shows
that in our situation, the 4--th moment is always bounded by a
constant times the square of the 2--nd moment. As was explained
above, this implies that the second moments remain bounded. By virtue
of \tht{8.3} this exactly means that $\Vert\om_n\Vert'$ remains
bounded, as was required. 
\enddemo

\demo{Step 4} Let a sequence $(\om_n)$ be such that $(F^{(\om_n)})$
is convergent. Here we will show that $\Vert\om_n\Vert$ remains
bounded. By virtue of step 2 this will imply that $(\om_n)$ converges
in $\om$. As we noticed in the very beginning of the section, our
topological spaces are separable and metrizable. Therefore, it will
follow that the inverse map $F^{(\om)}\mapsto\om$ is continuous, which
will complete the proof of the theorem.

By virtue of step 3, $\Vert\om_n\Vert'$ remains bounded. We would like
to deduce from this that $\Vert\om_n\Vert$ remains bounded, too. 

Let us compare $\Vert\om\Vert$ and $\Vert\om\Vert'$. We have
$$
\al^\pm_i\le\al^\pm_i(1+\al^\pm_i)
$$
and
$$
\be^\pm_i\le2\be^\pm_i(1-\be^\pm_i), \quad\text{provided that
$\be^\pm_i\in[0,\tfrac12]$.} 
$$
In general, we only know that the coordinates $\be^\pm_i$ are in
$[0,1]$, so that there is no universal estimate of the form
$\Vert\om\Vert\le\const\cdot\Vert\om\Vert'$. We will bypass this
difficulty as follows.

For any $\om\in\Om$, only finitely many coordinates $\be^+_i$ or
$\be^-_i$ are strictly greater than $\tfrac12$. Let us call these
coordinates ``bad'', and let $K(\om)$ denote their number. Then we have
an estimate of the form $\Vert\om\Vert\le C\cdot\Vert\om\Vert'$,
where $C$ depends only on $K(\om)$. Now it suffices to show that
$K(\om_n)$ remains bounded, which will imply that $\Vert\om_n\Vert$
remains bounded. 

Recall that for any $\om\in\Om$, $\be^+_1+\be^-_1\le1$. Hence, both
$\be^+_1$ and $\be^-_1$ cannot be ``bad''. Therefore, if
$K(\om)>0$, then the ``bad'' coordinates are either the first $K(\om)$
coordinates of $\be^+$ or the first $K(\om)$ coordinates of $\be^-$. 
Define a new element $\widetilde\om\in\Om$ as follows. If the ``bad''
coordinates were in $\be^+$ then we remove them, then add to $\be^-$ new
coordinates $1-\be^+_1,\dots,1-\be^+_{K(\om)}$, and finally rearrange
all coordinates in $\be^-$ in the descending order. If the ``bad''
coordinates were in $\be^-$ then we do the same operation with
$\be^+$ and $\be^-$ interchanged. Then we have  
$$
K(\widetilde\om)=0, \qquad \Vert\widetilde\om\Vert'=\Vert\om\Vert'.
\tag8.5
$$
On the other hand (see Remark 1.5),

$$
F^{(\widetilde\om)}(u)=F^{(\om)}(u)u^{\mp K(\om)}\,. \tag8.6
$$

Let us assume now that in our sequence $(\om_n)$, the quantity
$K(\om_n)$ is unbounded. Taking a subsequence, we may assume that
$K(\om_n)\to\infty$. Consider the corresponding sequence
$(\widetilde{\om_n})$. Since $\Vert\om_n\Vert'$ remains bounded and
$\Vert\widetilde{\om_n}\Vert'=\Vert\om_n\Vert'$ (see \tht{8.5}), we
conclude that $\Vert\widetilde{\om_n}\Vert'$ is bounded, too. Since
$\widetilde{\om_n}$ does not have ``bad'' coordinates (see \tht{8.5}),
the argument above shows that $\Vert\widetilde{\om_n}\Vert$ remains
bounded. It follows that the 
corresponding functions on the unit circle form a relatively compact
set in the topology of uniform convergence. Taking a subsequence
again, we may assume that the functions $F^{(\widetilde{\om_n})}$
converge. On the other hand, by the initial assumption, the functions
$F^{(\om_n)}$ converge, too. Now taking account of \tht{8.6} we
obtain a contradiction. Indeed, we get two sequences of continuous
functions on the unit circle, normalized at 1, say $F_n(u)$ and
$\widetilde F_n(u)$, such that
$$
\widetilde F_n(u)=F_n(u)u^{\mp K_n}\,, \qquad \text{where $K_n\to\infty$}, 
$$
and such that both $(F_n)$ and $(\widetilde F_n)$ converge in the
uniform metric. These conditions imply that the ratios
$$
\frac{\widetilde F_n(u)}{F_n(u)}=u^{\mp K_n}
$$
uniformly converge in a neighborhood of $u=1$, which is impossible,
because $K_n\to\infty$.  

This contradiction shows that actually $K(\om_n)$ must be bounded,
which completes the proof. \qed
\enddemo  
\enddemo

\head 9. Existence and uniqueness of spectral decomposition \endhead

The aim of this section is to rederive (and slightly refine) a result
due to Voiculescu:

\proclaim{Theorem 9.1 (Cf. \cite{Vo, Th\'eor\`eme 2})} For any
character $\chi$ of the group $U(\infty)$ there exists a probability
measure $P$ on the topological space $\Om$ such that
$$
\chi(U)=\int_{\Om}\chi^{(\om)}(U)P(d\om), \qquad U\in U(\infty).
$$
Moreover, such a measure is unique.
\endproclaim

We call $P$ the {\it spectral measure\/} of $\chi$. 

\demo{Comments} 1. Recall (see the beginning of \S8) that $\Om$ is a
good topological space (locally compact, separable, metrizable), so
that there is no problem in defining measures on it. Specifically, we
take Borel measures with respect to the natural Borel structure of
$\Om$.  

2. The integral above makes sense, because $\chi^{(\om)}(U)$ is a
continuous function in $\om$ (see Theorem 8.1) and
$|\chi^{(\om)}(U)|\le1$. 

3. It is readily seen that for any probability Borel measure $P$ on
$\Om$, the above formula defines a character, so that the
correspondence $\chi\mapsto P$ is a bijection between the space $\Cal
X(U(\infty))$ of characters and the space of probability Borel
measures on $\Om$. 

4. In Remark 9.4 below we compare our approach with that of
Voiculescu.  
\enddemo

We shall derive Theorem 9.1 from a more general claim, see Theorem
9.2 below. 

Assume we are given a sequence $\Ga_0,\Ga_1,\Ga_2,\dots$ of nonempty
sets, where $\Ga_0$ consists of a single point denoted by the symbol
$\varnothing$ while each $\Ga_N$ with $N\ge1$ is a finite or
countable set. Let $\De_N$ denote the space of formal convex
combinations of the points of $\Ga_N$; this is a simplex whose
vertices are points of $\Ga_N$. Further, assume we are given a
function $q(\nu,\la)$ defined on couples
$(\nu,\la)\in\Ga_{N-1}\times\Ga_N$, where $N=1,2,\dots$, such that
$0\le q(\nu,\la) \le1$ and, for any $\la\in\Ga_N$, $\sum_\nu
q(\nu,\la)=1$. In the particular case $N=1$ this means that
$q(\varnothing,\la)=1$ for all $\la\in\Ga_1$.

For each $N=1,2,\dots$, there exists a unique affine map
$\De_N\to\De_{N-1}$ taking any $\la\in\Ga_N$ to the convex
combination $\sum_\nu q(\nu,\la)\nu\in\De_{N-1}$. Let
$\De=\varprojlim\De_N$ denote 
the projective limit taken with respect to these affine maps. 

Let $\Ga=\Ga_0\sqcup\Ga_1\sqcup\Ga_2\sqcup\dots$ be the disjoint
union of the sets $\Ga_N$; this is a countable set. Denote by $\F$ 
the vector space of all real functions on $\Ga$ and equip it with the
topology of pointwise convergence on $\Ga$. Note that $\F$ is a locally
convex vector space. It is metrizable, because $\Ga$ is countable.

We identify $\De$ with the
subset of $\F$ formed by nonnegative functions $f$ such that
$f(\varnothing)=1$ and $f(\nu)=\sum_{\la}q(\nu,\la)f(\la)$ for any
$\nu\in\Ga_{N-1}$, $N=1,2,\dots$, where the summation is taken over
$\la\in\Ga_N$ (note that these conditions imply that 
$\sum_\la f(\la)=1$, summed over $\la\in\Ga_N$). 

Clearly, $\De$ is a convex subset of $\F$. As a subset of $\F$, $\De$
inherits its topology. We also consider the Borel structure on $\De$
generated by this topology. One can prove that this structure is standard.

Let the symbol $\Ex(\,\cdot\,)$ denote the subset of extreme points
of a convex set.

\proclaim{Theorem 9.2} Let $\De$ be the convex set $\varprojlim\De_N$
defined above. The subset $\Ex(\De)$ of extreme points of
$\De$ is a Borel subset and each point $f\in\De$ is uniquely
represented by a probability Borel measure $P$ on $\Ex(\De)$:
$f=\int_{\Ex(\De)}gP(dg)$. That is,
$f(\la)=\int_{\Ex(\De)}g(\la)P(dg)$ for any $\la\in\Ga$.
\endproclaim

\demo{Derivation of Theorem 9.1 from Theorem 9.2} Take $\Ga_N=\GT_N$,
$N=1,2,\dots$ and take as $q(\nu,\la)$ the cotransition probability
function (Definition 7.2). Then $\De_N$ turns into the set of probability
measures on $\GT_N$ and the set $\De$ becomes the set of coherent
families of measures on $\GT$ (Definition 7.3). By virtue of
Proposition 7.4 we get a bijection between $\De$ and $\Cal X(U(\infty))$,
which is an isomorphism of convex sets. Moreover, this is a
homeomorphism of topological spaces (see the discussion of
the topology on $\Cal X(U(\infty))$ in the beginning of \S8).
Applying Theorem 1.3, which 
provides an explicit parameterization of extreme characters, we get a
bijective correspondence between $\Ex(\De)$ and $\Om$. By Theorem
8.1, this correspondence is a homeomorphism of topological spaces.
Hence, it preserves the Borel structures. This turns Theorem 9.1 into
a special case of Theorem 9.2. \qed
\enddemo

The proof of Theorem 9.2 is based on the following lemma.

\proclaim{Lemma 9.3} The convex set $\De$ is a Choquet simplex, i.e.,
the cone generated by $\De$ is a lattice.
\endproclaim

\demo{Proof} Let $C$ denote this cone. It coincides with
the subset of $\F$ that is described similarly to $\De\subset\F$: the
only difference is that the condition $f(\varnothing)=1$ is dropped. 

We extend the function $q(\nu,\la)$ to any couples $\nu\in\Ga_M$,
$\la\in\Ga_N$ with $M<N$ using the following recurrence relation:
$$
q(\nu,\la)=\sum_{\mu\in\Ga_{N-1}}q(\nu,\mu)q(\mu,\la).
$$

Given $f_1,f_2\in C$ we construct their lowest upper bound as follows.
Define a function $f$ on $\Ga$ by
$$
f(\nu)=\lim_{N\to\infty}\sum_{\la\in\Ga_N}q(\nu,\la)\max(f_1(\la),f_2(\la)).
$$
The limit exists, because, for any fixed $\nu$, the $N$th sum 
monotonically increases as $N\to\infty$ and remains bounded from
above by $f_1(\nu)+f_2(\nu)$. It is readily verified that $f$ belongs
to the cone and is the lowest upper bound for $f_1$ and $f_2$.

The existence of the greatest lower bound is verified similarly: it
suffices to substitute ``$\min$'' for ``$\max$''. \qed
\enddemo

\demo{Proof of Theorem 9.2} If all sets $\Ga_N$ are finite, then
$\De$ is compact. Since $\De$ is metrizable (as $\F$ is metrizable),
the claims of Theorem 9.2 immediately follow from Lemma 9.3 and
Choquet's theorem, see \cite{Ph}. When 
the sets $\Ga_N$ are allowed to be countable, the space $\De$ may be
noncompact, so that Choquet's theorem is not directly applicable. To
overcome this difficulty we embed $\De$ into a bigger set
$\wt\De\subset\F$ which is compact. 

Specifically, let $\wt\De$ be the set of nonnegative functions $f$ on
$\Ga$ such that $f(\varnothing)=1$ and
$f(\nu)\ge\sum_{\la}q(\nu,\la)f(\la)$ for any $\nu\in\Ga_{N-1}$,
$N=1,2,\dots$, where $\la$ ranges over $\Ga_N$ (i.e., in the
definition of $\De$, we have replaced the equality by an inequality). We
note that $\wt\De$ is a {\it compact\/} convex set containing $\De$.

Next, for any $M=0,1,\dots$, let $\wt\De_M$ be the set of nonnegative
functions on $\Ga$ such that $f(\varnothing)=1$,
$f(\nu)=\sum_{\la}q(\nu,\la)f(\la)$ for all $\nu\in\Ga_{N-1}$ with
$N\le M$, and $f(\,\cdot\,)\equiv0$ on
$\Ga_{M+1}\sqcup\Ga_{M+2}\sqcup\dots$. Note that $\wt\De_M$ is a
convex subset of $\wt\De$, isomorphic to $\De_M$. 

Finally, dropping the condition $f(\varnothing)=1$ in the definitions
above we get the cones spanned by the sets $\wt\De$ and $\wt\De_M$;
we denote them by $\wt C$ and $\wt C_M$, respectively.  

Now the crucial remark is that any element $f\in\wt C$ is uniquely
represented as a sum $f=f_0+f_1+\dots+f_\infty$, where 
$f_M\in C_M$ and $f_\infty\in C$. Conversely, any such sum
represents an element of $\wt C$ provided that 
$\sum_M f_M(\varnothing)<+\infty$, cf. \cite{Ol3, \S22}.

This implies a number of consequences. First, the cone $\wt C$ is a
lattice, so that we may apply Choquet's theorem to $\wt\De$. Next,
$\Ex(\wt\De)$ is the disjoint union of the sets
$\Ex(\wt\De_M)\simeq\Ga_M$, where $M=0,1,\dots$, and the set
$\Ex(\De)$. Thus, $\Ex(\De)$ is the difference of $\Ex(\wt\De)$,
which is a Borel set, and a countable set; whence, $\Ex(\De)$ is a
Borel set. Finally, if $f$ belongs to $\De$ then its representing
measure on $\Ex(\wt\De)$ is concentrated on
$\Ex(\De)\subset\Ex(\wt\De)$. This concludes the proof.  \qed
\enddemo

\example{Remark 9.4} Here is a comment to Voiculescu's original
result \cite{Vo, Th\'eor\`eme 2}. Its formulation says that any
$\chi\in\Cal X(U(\infty))$ is uniquely represented by a 
probability measure $P$ on the space $\Ex(\Cal X(U(\infty)))$, equipped
with the Borel structure inherited from the ambient space $\Cal
X(U(\infty))$. To prove the existence of $P$, Voiculescu also embedded
$\Cal X(U(\infty))$ into a compact space. But he used a special property
of the cotransition function $q(\nu,\la)$ of the graph $\GT$
(specifically, the fact that $q(\nu,\la )$ tends to zero
when $\nu\in\GT_{N-1}$ is fixed and $\la\in\GT_N$ goes to infinity).
This allowed him to take, instead of our space $\widetilde\De$, a
smaller set. Voiculescu's proof of the uniqueness statement is quite
different from ours: it substantially relies on the multiplicativity
property of extreme characters. The argument presented above seems to
be more general and direct.  

Voiculescu's original formulation did not involve the space $\Om$,
because at that time it was not yet clear whether the characters
$\chi^{(\om)}$ exhaust the whole set $\Ex(\Cal X(U(\infty)))$ of
extreme characters. Our (modest) supplement consists in checking that
the bijection between the spaces $\Ex(\Cal X(U(\infty)))$ and $\Om$
preserves the Borel structures, so that $P$ can be carried over from
the ``abstract'' space $\Ex(\Cal X(U(\infty)))$ to the ``concrete''
space $\Om$. Of course, this is a pure technical claim whose validity
seems to be beyond doubt. However, it is not completely trivial.  
\endexample

\head 10. Approximation of spectral measures \endhead

The aim of this section is to establish a relationship between
the spectral measure $P$ of an arbitrary character $\chi$ and the
coherent system $\{P_N\}$ corresponding to $\chi$. Recall that $P_N$
is a probability measure on the discrete set $\GT_N\subset\Z^N$
($N=1,2,\dots$) while $P$ is a probability Borel measure on the region
$\Om\subset\R^{4\infty+2}$ (see Proposition 7.4 and Theorem 9.1). We
will prove that, as $N\to\infty$, the measures $P_N$ approximate the
measure $P$ in a certain sense. To precisely state this claim we need
a few definitions.  

Given a signature $\la\in\GT_N$, we denote by $\la^+$ and $\la^-$ its
positive and negative parts. These are two Young diagrams such that
$\ell(\la^+)+\ell(\la^-)\le N$, where $\ell(\,\cdot\,)$ is the number
of nonzero rows of a Young diagram. That is,
$$
\la=(\la^+,\dots,\la^+_k,0,\dots,0,-\la^-_l,\dots,-\la^-_1), \qquad
k=\ell(\la^+), \quad l=\ell(\la^-).
$$

Next, given a Young diagram $\nu$, we denote by $d(\nu)$ the number
of diagonal boxes in $\nu$, and we introduce the {\it Frobenius
coordinates\/} of $\nu$:
$$
p_i(\nu)=\nu_i-i, \quad q_i(\nu)=\nu'_i-i, \qquad i=1,\dots,d(\nu),
$$
where $\nu'$ stands for the transposed (conjugate) diagram. 

Then, following Vershik--Kerov \cite{VK1}, we introduce the {\it
modified Frobenius coordinates\/} of $\nu$ as follows
$$
\wt p_i(\nu)=p_i(\nu)+\tfrac12=\nu_i-i+\tfrac12, \quad 
\wt q_i(\nu)=q_i(\nu)+\tfrac12=\nu'_i-i+\tfrac12, 
\qquad i=1,\dots,d(\nu).
$$
Note that
$$\gather
\wt p_1(\nu)>\dots>\wt p_{d(\nu)}(\nu)>0, \quad
\wt q_1(\nu)>\dots>\wt q_{d(\nu)}(\nu)>0, \\
\sum_{i=1}^{d(\nu)}(\wt p_i(\nu)+\wt q_i(\nu))=|\nu|,
\endgather
$$
where $|\nu|$ denotes the total number of boxes in $\nu$. We also
agree to set
$$
\wt p_i(\nu)=\wt q_i(\nu)=0, \quad i=d(\nu)+1, \, d(\nu)+2, \dots,
$$
so that the coordinates $\wt p_i(\nu)$ and $\wt q_i(\nu)$ are now
defined for any $i=1,2,\dots$\,.

\example{Definition 10.1} For any $N=1,2,\dots$, we embed the set
$\GT_N$ into $\Om$ as follows
$$
\gather
\GT_N\ni\la\,\longmapsto\,
\om=(a^+,b^+,a^-,b^-,c^+,c^-)\in\Om,  \\
a^\pm_i=\frac{\wt p_i(\la^\pm_N)}N\,, \quad
b^\pm_i=\frac{\wt q_i(\la^\pm_N)}N \quad (i=1,2,\dots), \quad
c^\pm=\frac{|\la^\pm_N|}N\,.
\endgather
$$
\endexample

\proclaim{Theorem 10.2} Let $\chi$ be an arbitrary character of
$U(\infty)$, $P$ be its spectral measure, and $\{P_N\}$ be the
coherent system corresponding to $\chi$. Further, for any $N$,
convert $P_N$ to a probability measure $\un P_N$ on $\Om$, which is
the pushforward of $P_N$ under the embedding $\GT_N\hookrightarrow\Om$
introduced in Definition 10.1.

Then, as $N\to\infty$, the measures $\un P_N$ weakly converge to the
spectral measure $P$. I.e., for any bounded continuous function $F$
on $\Om$,
$$
\lim_{N\to\infty}\int_\Om F(\om)\un P_N(d\om)=
\int_\Om F(\om) P(d\om).
$$
\endproclaim 

The proof is quite similar to that of \cite{BO3, Theorem 5.3}. 
We will see that Theorem 10.2 is a corollary of another general
result, Theorem 10.7.

A {\it path\/} in the Gelfand--Tsetlin graph $\GT$ 
is an infinite sequence $t=(t_1,t_2,\dots)$ such that $t_N\in\GT_N$ and
$t_N\prec t_{N+1}$ for any $N=1,2,\dots$. The set of the paths will
be denoted by $\Cal T$. 

We also need {\it finite paths.\/} A finite path of length $N$ is a
sequence $\tau=(\tau_1,\dots,\tau_N)$, where $\tau_1\in\GT_1$,\dots,
$\tau_N\in\GT_N$ and $\tau_1\prec\dots\prec\tau_N$. The set of
finite paths of length $N$ will be denoted by $\Cal T_N$. One can
identify $\Cal T$ with the projective limit space 
$\varprojlim\Cal T_N$, where the projection $\Cal T_N\to\Cal T_{N-1}$
is the operation of removing the last vertex $\tau_N$. 

Consider the natural embedding $\Cal T\subset\prod\limits_N \GT_N$.
We equip $\prod\limits_N \GT_N$ with the product topology (the sets
$\GT_N$ are viewed as discrete spaces). The set $\Cal T$ is
closed in this product space. We equip $\Cal T$ with the induced
topology. Equivalently, the topology is that of the projective limit
space $\varprojlim\Cal T_N$. Then $\Cal T$ turns into a totally
disconnected topological space.  

Given a finite path $\tau=(\tau_1,\dots,\tau_N)\in\Cal T_N$, define
the {\it cylinder set\/} $C_\tau\subset\Cal T$ as the inverse image
of $\{\tau\}$ under the projection $\Cal T\to\Cal T_N$\,,
$$
C_\tau=\{t\in\Cal T\mid t_1=\tau_1,\dots,t_N=\tau_N\}.
$$
The cylinder sets form a base of topology in $\Cal T$. 

Consider an arbitrary signature $\la\in\GT_N$. The set of finite paths 
$\tau=(\tau_1\prec\dots\prec\tau_N)$ ending at $\la$ has cardinality
equal to $\Dim_N\la=\chi^\la(e)$. The cylinder sets $C_\tau$
corresponding to these finite paths $\tau$ are pairwise disjoint, and
their union coincides with the set of infinite paths $t$ passing
through $\la$.  

A {\it central measure\/} is any probability Borel measure on $\Cal
T$ such that the mass of any cylinder set $C_\tau$ depends only of 
its endpoint $\la$. Clearly, central measures form a convex set. 

These definitions are inspired by \cite{VK1}. 

\proclaim{Proposition 10.3} There exists a natural bijective
correspondence $M\longleftrightarrow\{P_N\}$ 
between central measures $M$ and coherent systems $\{P_N\}$, defined
by the relations 
$$
\Dim_N\la\cdot M(C_\tau)=P_N(\la), \tag10.1
$$
where $N=1,2,\dots$, $\la\in\GT_N$, and $\tau$ is an arbitrary finite
path ending at $\la$. 
\endproclaim

In other words, the relations mean that for any $N$, the
pushforward of $M$ under the natural projection 
$$
\prod_{N=1}^\infty\GT_N\supset \Cal T\to \GT_N \tag10.2
$$
coincides with $P_N$.

\demo{Proof} Let $\{P_N\}$ be a coherent system. For any $N$, we
define a measure $M_N$ on the discrete space $\Cal T_N$
as follows. Given $\tau\in\Cal T_N$, we set
$$
M_N(\{\tau\})=\frac1{\Dim_N\la}\,P_N(\la),
$$
where $\la\in\GT_N$ is the end of $\tau$. This is a probability
measure. Its pushforward under the projection $\Cal T_N\to\Cal
T_{N-1}$ coincides with $M_{N-1}$: indeed, this is exactly a
reformulation of the coherency property \tht{7.2}. Hence the family
$\{M_N\}$ determines a probability measure $M$ on the projective
limit space $\Cal T$. Clearly, $M$ is a central measure, and we have
the relations \tht{10.1}.

Conversely, let $M$ be a central measure. For any $N$, define a
probability measure $P_N$ on $\GT_N$ as the pushforward of $M$ under
the projection \tht{10.2}. The fact that $M$ is central then implies
that the family $\{P_N\}$ satisfies the coherency property.  \qed
\enddemo

\proclaim{Corollary 10.4} There is a natural bijective correspondence
$\chi\longleftrightarrow M$ between characters and central measures.
This correspondence is an isomorphism of convex sets.
\endproclaim

\demo{Proof} The bijection is established by means of the bijections
$\chi\longleftrightarrow\{P_N\}$ (Proposition 7.4) and 
$\{P_N\}\longleftrightarrow M$ (Proposition 10.3). {}From the proofs
of these propositions it is clear that this is an isomorphism of
convex sets. \qed
\enddemo

\example{Definition 10.5} Let $t=(t_N)\in\Cal T$ be a path. We say that
$t$ is {\it regular\/} if the images of the $t_N$'s under the
emeddings $\GT_N\hookrightarrow\Om$ (see Definition 10.1) converge to
a point $\om\in\Om$. Then $\om$ is called the {\it end\/} of $t$.
Equivalently, the condition means that there exist limits
$$
\al^\pm_i=\lim_{N\to\infty}\frac{\wt p_i(t^\pm_N)}N\,, \quad
\be^\pm_i=\lim_{N\to\infty}\frac{\wt q_i(t^\pm_N)}N \quad (i=1,2,\dots), \quad
\de^\pm=\lim_{N\to\infty}\frac{|t^\pm_N|}N\,,
$$
where $t^+_N$ and $t^-_N$ denote the positive and negative parts of
$t_N\in\GT_N$. The set of regular paths will be denoted by $\Treg$.
\endexample

Introduce the map
$$
\pi:\Treg\to\Om, \qquad 
\Treg\ni t\longmapsto\{\text{the end of $t$}\}\in\Om. \tag10.3
$$
Further, for any $N=1,2,\dots$, we introduce the map $\pi_N:\Cal
T\to\Om$ as the composition 
$$
\pi_N:\Cal T@>{t\mapsto t_N}>>\GT_N@>{\text{Definition 10.1}}>>\Om.
\tag10.4 
$$
For any $t\in\Treg$, $\pi_N(t)$ converges (as $N\to\infty$) to 
$\pi(t)$ in the topology of the space $\Om$. Indeed, this holds by
the very definition of regular paths, see Definition 10.5. 

\proclaim{Lemma 10.6} $\Treg\subset\Cal T$ is a Borel set, and the
maps $\pi$, $\pi_N$ are Borel maps from $\Treg$ to $\Om$.
\endproclaim

\demo{Proof} Indeed, for any $N$, the functions
$$
t\mapsto\frac{\wt p_i(t^\pm_N)}N\,, \quad
t\mapsto\frac{\wt q_i(t^\pm_N)}N \quad (i=1,2,\dots), \quad
t\mapsto\frac{|t^\pm_N|}N
$$
are continuous functions on the space $\Cal T$. This implies that
the regularity condition determines a set of type $F_{\sigma\de}$,
hence a Borel set. 

Each $\pi_N$, being a cylindrical map, is continuous on $\Cal T$.
Its restriction to $\Treg$ is also continuous, hence Borel. Finally,
$\pi$ is the pointwise limit of the $\pi_N$'s, hence it is a Borel
map, too. 
\enddemo

\proclaim{Theorem 10.7} Any central measure $M$ is concentrated on the
set $\Treg\subset\Cal T$. The pushforward $\pi(M)$ under the
projection \tht{10.3} coincides with the spectral measure $P$ of the
character $\chi$ that corresponds to $M$. 
\endproclaim 

This claim makes sense, because $\Treg$ is a Borel set.

\demo{Derivation of Theorem 10.2 from Theorem 10.7} Let $M$ be the
central measure corresponding to $\chi$.

By Proposition 10.3, we have $\pi_N(M)=\un P_N$. Now we will interpret
$M$ as a probability measure on $\Treg$, which is possible by virtue
of Theorem 10.7. Also, $\pi_N$ will be viewed as a projection
$\Treg\to\Om$. Then the equality $\pi_N(M)=\un P_N$ remains true. On the
other hand, $\pi(M)=P$ (again by Theorem 10.7). Thus, both $\un P_N$ and
$P$ can be viewed as pushforwards of one and the same probability
measure $M$. 

Now let $F$ be a bounded continuous function on $\Om$. Since
$\pi_N(t)\to\pi(t)$ for any $t\in\Treg$, 
$$
F(\pi_N(t))\to F(\pi(t)), \qquad t\in\Treg.
$$
Thus, we get a sequence $\{F(\pi_N(\,\cdot\,))\}$ of uniformly bounded
functions on $\Treg$ converging pointwise to the function
$F(\pi(\,\cdot\,))$. All these functions are Borel functions, because
$\pi_N$ and $\pi$ are Borel maps (Lemma 10.6).  Consequently, as
$N\to\infty$, 
$$
\int_{\Treg}F(\pi_N(t))M(dt)\, \longrightarrow \,
\int_{\Treg} F(\pi(t)) M(dt).
$$
Since $\pi_N(M)=\un P_N$ and $\pi(M)=P$, we can convert these
integrals to integrals over $\Om$,
$$
\int_{\Treg}F(\pi_N(t))M(dt)=\int_\Om F(\om)\un P_N(d\om), \qquad
\int_{\Treg} F(\pi(t)) M(dt)=\int_\Om F(\om) P(d\om),
$$
which concludes the proof. \qed
\enddemo

The rest of the section is devoted to the proof of Theorem 10.7. 

Given two signatures, $\nu\in\GT_{n}$ and $\la\in\GT_N$, where
$n<N$, denote by $\Dim_{nN}(\nu,\la)$ the number of paths
starting at $\nu$ and ending at $\la$. I.e., the number of chains
$$
\tau=(\tau_{n}\prec\tau_{n+1}\prec\dots
\prec\tau_{N-1}\prec\tau_N), \qquad \tau_{n}=\nu, \quad \tau_N=\la.
$$
We extend the definition of the cotransition function by setting (cf.
the proof of Lemma 9.3)
$$
q(\nu,\la)=\frac{\Dim_n(\nu)\Dim_{nN}(\nu,\la)}{\Dim_N\la}\,,
\qquad \nu\in\GT_{n}, \quad \la\in\GT_N, \quad n<N.
$$ 
Note that
$$
q(\nu,\la)=\sum_\tau \prod_{i=n+1}^N q(\tau_{i-1},\tau_i),
$$
summed over all $\tau$'s as above. 

For any coherent system $\{P_N\}$, we have
$$
P_{n}(\nu)=\sum_{\la\in\GT_N}q(\nu,\la)P_N(\la), 
\qquad \nu\in\GT_{n}, \quad n<N. 
$$
Indeed, this relation is obtained by iterating the coherency relation
\tht{7.2}. 

\proclaim{Proposition 10.8} Let $\{P_N\}$ be a coherent system and
$M$ be the corresponding central measure. Assume that $M$ is extreme.
Then for $M$--almost all paths $t=(t_N)\in\Cal T$
$$
\lim_{N\to\infty}q(\nu,t_N)=P_{n}(\nu)\,, \qquad 
n=1,2,\dots, \quad \nu\in\GT_{n}\,. \tag10.5
$$
\endproclaim

\demo{Comment} This result actually holds for general ``branching
graphs'' in the sense of Vershik--Kerov \cite{VK3}. It was stated
in their paper \cite{VK1} for the Young graph, associated with the
infinite symmetric group $S(\infty)$. A detailed proof (an adaptation
of the proof of the Birkhoff--Khinchine ergodic theorem) is contained in
an unpublished work by Kerov. We present below a similar argument but
use a limit theorem for reversed martingales (as was suggested in
\cite{VK1}). Note also a very close earlier result due to Vershik, see
\cite{Ve, Theorem 1} and  \cite{OV, Theorem 3.2 and Remark 3.6}.
\enddemo

\demo{Proof} {\it Step 1.\/} Let us say that two paths
$t=(t_m)$, $t'=(t'_m)$ are {\it $N$--equivalent\/} if $t_m=t'_m$ for
$m\ge N$. Denote by $\xi_N(t)$ the $N$--equivalence class containing
$t$. Then
$$
\{t\}=\xi_1(t)\subset\xi_2(t)\subset\dots\,.
$$
We say that $t$ and $t'$ are {\it $\infty$--equivalent\/} if they are
$N$--equivalent for $N$ large enough. The $\infty$--equivalence class
of $t$ is denoted by $\xi_\infty(t)$. Clearly,
$\xi_\infty(t)=\cup\xi_N(t)$.  

Let $\Cal B_{-N}$ denote the $\sigma$--algebra of those Borel sets in
$\Cal T$ that are saturated with respect to the $N$--equivalence
relation, where $N=1,2,\dots$\,. We define $\Cal B_{-\infty}$
likewise. We have 
$$
\Cal B_{-\infty}\subset\dots\subset\Cal B_{-2}\subset\Cal B_{-1}
=\{\text{All Borel sets}\}, 
\qquad \Cal B_{-\infty}=\bigcap\limits_{N=1}^\infty\Cal B_{-N}
$$

Fix an arbitrary bounded Borel function $\psi$ on $\Cal T$. We will
view $\psi$ as a random variable defined on the probability space
$(\Cal T, M)$. Let $\psi_N=\E(\psi\mid\Cal B_{-N})$ and 
$\psi_\infty=\E(\psi\mid\Cal B_{-\infty})$ be the conditional
expectations of $\psi$ with respect to the $\sigma$--algebras 
$\Cal B_{-N}$ and $\Cal B_{-\infty}$, respectively. 

By the theorem on convergence of (reversed) martingales, we have
$$
\lim_{N\to\infty}\psi_N=\psi_\infty \qquad \text{almost everywhere
with respect to $M$,}
$$
see, e.g., Doob \cite{Do, Chapter VII, Theorem 4.2}. 

{\it Step 2.\/} So far we did not use the assumption that $M$ is
central. Now we remark that this assumption makes it possible to
describe each $\psi_N$ explicitly as the averaging over the
$N$--equivalence classes:
$$
\psi_N(t)=\frac1{\Dim_N(t_N)}\sum_{t'\in\xi_N(t)}\psi(t'). \tag10.6
$$

Next, we use the assumption that $M$ is an extreme central measure to
conclude that each set $A\in\Cal B_{-\infty}$ has mass 0 or 1.
Indeed, assume $0<M(A)<1$, and let $\boldkey{1}_A$ stand for the
characteristic function of $A$. Then the measure $M$ can be written as
a nontrivial convex combination of two probability measures,
$$
M=M(A)\cdot\frac{\boldkey{1}_A\,M}{M(A)}+
(1-M(A))\frac{\boldkey{1}_{\Cal T\setminus A}\,M}{1-M(A)}\,. 
$$
Since $A$ is saturated with respect to the $\infty$--equivalence
relation, these two measures are central, which contradicts the
extremality assumption. 

It follows that $\psi_\infty$ is constant almost everywhere. Since the
means of all $\psi_N$'s are the same and equal to the mean of $\psi$,
the constant is the mean of $\psi$. 

Therefore, we conclude: let $\psi$ be any bounded Borel function on
$\Cal T$ and let $\psi_N(t)$ be defined by \tht{10.6}; then
$$
\lim_{N\to\infty}\psi_N(t)=\int_{\Cal T}\psi(t')M(dt') \qquad
\text{for $M$--almost all $t\in\Cal T$}. \tag10.7
$$

{\it Step 3.\/} Fix $n$ and $\nu\in\GT_n$, and set
$$
\psi(t)=\cases 1, & \text{if $t_n=\nu$}, \\ 0, & \text{otherwise.} 
\endcases
$$
For any $N\ge n$, we have
$$
\psi_N(t)=\frac{\Dim_n(\nu)\Dim_{nN}(\nu, t_N)}{\Dim_N t_N}
=q(\nu,t_N).
$$
Indeed, the set $\xi_N(t)$ contains exactly $\Dim_n(\nu)\Dim_{nN}(\nu, t_N)$
paths $t'$ passing through $\nu$.
Finally, the mean of $\psi$ equals $P_n(\nu)$. 

Now we apply the result of Step 2 and get from \tht{10.7} the required
claim. \qed 
\enddemo

Assume that for any $N=1,2,\dots$ we are given a character (in the
sense of Definition 1.1) of the group $U(N)$, denoted as $\chi^\n$.
Let us say that the sequence $\{\chi^\n\}$ {\it converges\/} to a
function $\chi^{(\infty)}$ defined on the group $U(\infty)$ if
$$
\text{as $N\to\infty$,}\quad \chi^\n\mid_{U(n)}\,
@>{\text{uniformly}}>>\,
\chi^{(\infty)}\mid_{U(n)} \,,
\quad \text{ for any fixed $n=1,2,\dots$} \tag10.8
$$
Clearly, the limit function $\chi^{(\infty)}$ is a character of
$U(\infty)$. Expand $\chi^\n$ on irreducible normalized characters of
$U(N)$:
$$
\chi^\n=\sum_{\la\in\GT_N}P^\n \wt\chi^{\,\la}.
$$
Then, for any $n<N$,
$$
\chi^\n\mid_{U(n)}=\sum_{\nu\in\GT_{n}}
\left(\sum_{\la\in\GT_N}q(\nu,\la)P^\n(\la)\right)\wt \chi^{\,\nu}
$$
It follows that \tht{10.8} is equivalent to the following condition: for
any $n$ and any $\nu\in\GT_{n}$,
$$
\lim_{N\to\infty}\sum_{\la\in\GT_N}q(\nu,\la)P^\n(\la)=P_{n}(\nu),
\tag10.9 
$$
where $\{P_N\}$ is the coherent system corresponding to $\chi^{(\infty)}$.

\proclaim{Proposition 10.9} Let $\{\la(N)\in\GT_N\}_{N=1,2,\dots}$ be
a sequence of signatures and let $\chi^\n=\wt\chi^{\,\la(N)}$ be the
corresponding normalized irreducible characters.

The sequence $\{\chi^\n\}$ converges to a function $\chi^{(\infty)}$
if and only if the images of the $\la(N)$'s under the embeddings of
Definition 10.1 converge to a point $\om\in\Om$, and then the limit
function $\chi^{(\infty)}$ coincides with the extreme character
$\chi^{(\om)}$.
\endproclaim

\demo{Proof} This result is due to Vershik--Kerov \cite{VK2}. For a
detailed proof, see \cite{OkOl}. \qed
\enddemo

\proclaim{Proposition 10.10} Let $\om\in\Om$ and $M^{(\om)}$ be the
extreme central measure corresponding to the extreme character
$\chi^{(\om)}$. The measure $M^{(\om)}$ is concentrated on the subset
$\pi^{-1}(\om)\subset\Treg$ of regular paths ending at $\om$. 
\endproclaim

\demo{Proof} By Proposition 10.8, the measure $M^{(\om)}$ is concentrated
on the set of paths $t=(t_N)\in\Cal T$ satisfying the condition
\tht{10.5}, where $\{P_N\}$ is the coherent system corresponding to the
character $\chi^{(\om)}$. Let us show that this set coincides with
$\pi^{-1}(\om)$. 

Indeed, the condition \tht{10.5} coincides with the condition
\tht{10.9} for the characters $\chi^\n=\wt\chi^{\,t_N}$ (here $P^\n$ is
reduced to the delta measure at $t_N$). As explained
above, the latter condition is equivalent to the convergence of the
characters $\chi^\n$ to the character $\chi^{(\om)}$. By Proposition 10.9,
this exactly means that $t$ is a regular path ending at $\om$. \qed
\enddemo

\demo{Proof of Theorem 10.7} Translating Theorem 9.1 into the language
of central measures we get the decomposition
$$
M=\int_{\Om}M^{(\om)} P(d\om).
$$ 
By Proposition 10.10, any extreme measure $M^{(\om)}$ is concentrated on
the subset $\pi^{-1}(\om)\subset\Treg$. Hence, $M$ is concentrated on
$\Treg$. This proves the first claim of Theorem 10.7. The second claim
also follows from the above decomposition and Proposition 10.10. \qed 
\enddemo

\head 11. Conclusion: the problem of harmonic analysis \endhead

Now we are in a position to state the problem of harmonic analysis on
the group $U(\infty)$:
\medskip

{\it Let $(\zw)\in\Dadm$ be an arbitrary admissible quadruple of
parameters (Proposition 7.10), let
$P_N=\{P_N(\,\cdot\,\mid\zw)\}_{N=1,2,\dots}$ be the coherent system
as defined in \tht{7.8}, let $\chi_\zw$ be the corresponding
character, and let $P$ be the spectral measure of $\chi_\zw$. 

Describe explicitly the measure $P$.}

\medskip

Theorem 10.2 suggests the idea to evaluate $P$ by means of the limit
transition from the measures $\un P_N$. This idea is realized in the
subsequent paper \cite{BO4}.   

Let $\chi$ be an arbitrary character of $U(\infty)$ and $T_\chi$ be
the corresponding spherical representation of $(G,K)$, see \S2. One
can prove that the spectral measure of $\chi$ also determines the
decomposition of $T_\chi$ into a (multiplicity free) direct integral
of irreducible spherical representations. Recall that in \S3, we
constructed a family $\{T_{zw}\}$ of unitary representations. If
$\Re(z+w)>-\tfrac12$ then $T_{zw}$ possesses a distinguished vector,
and if, moreover, both $z$ and $w$ are nonintegral, then this vector
is cyclic (Proposition 6.9), so that $T_{zw}$ coincides with the
representation $T_\chi$ with the character $\chi=\chi_{zw}$. Thus, in
this case, the spectral measure of the 
character $\chi_{zw}$ also governs the decomposition of $T_{zw}$.

If $\Re(z+w)>-\tfrac12$ but at least one of the parameters
$z,w$ is integral, then the spectral measure of the character
$\chi_{zw}$ refers to a proper subrepresentation of $T_{zw}$. To
decompose the whole representation $T_{zw}$ we need additional tools.
Cf. \cite{KOV}.

Finally, recall that the construction of the representations $T_{zw}$
makes sense even when $\Re(z+w)\le-\tfrac12$, although then the
characters $\chi_{zw}$ disappear. It would be interesting to study the
decomposition problem in this case as well. Cf. a similar problem
concerning infinite measures $m^{(s)}$ stated in \cite{BO3, \S8}.

\enddocument
\bye